\documentclass[12pt,a4paper,reqno]{amsart}

\usepackage{geometry}
\usepackage{mathtools}  
\usepackage{url}
\usepackage[hidelinks,breaklinks, urlcolor=black]{hyperref}
 
\usepackage[myheadings]{fullpage}  
\usepackage[english]{babel} 
\usepackage{enumerate}   

\usepackage{bm} 
\usepackage{mathrsfs}

\usepackage[full]{textcomp} 
\usepackage{newtxtext} 
\usepackage{cabin} 
\usepackage{zlmtt}
\usepackage[bigdelims,vvarbb]{newtxmath} 
\usepackage[cal=boondoxo]{mathalfa} 
 
\usepackage{microtype} %
  
\newcommand{\C}{\mathbb{C}}

\newcommand{\N}{\mathbb{N}}
\newcommand{\R}{\mathbb{R}} 

\newcommand{\Z}{\mathbb{Z}}

\newcommand{\Odi}[1]{\Odip{}{#1}}

\newcommand{\Odip}[2]{\mathcal{O}_{#1}\left(#2\right)}

\newcommand{\odip}[2]{{o}_{#1}\left(#2\right)}
\newcommand{\odi}[1]{\odip{}{#1}}

\newcommand{\ULI}{\text{ULI}}
\newcommand{\LLI}{\text{LLI}}

\allowdisplaybreaks

\renewcommand{\qedsymbol}{$\square$}
\newenvironment{Proof}[1][Proof]{\par\noindent\textbf{#1.}~}
{\hfill\qedsymbol\smallskip\par}

\newtheoremstyle{sltheorems}
{10pt}
{6pt}
{\slshape}
{}
{\bfseries}
{.}
{.5em}
{\thmname{#1}\thmnumber{ #2}\thmnote{ (#3)}}

\theoremstyle{sltheorems} 
\newtheorem{Theorem}{Theorem}
  
\newtheorem{Proposition}{Proposition}  

\newtheoremstyle{remark}
{10pt}
{6pt}
{\rm} 
{}
{\bfseries}
{.}
{.5em}
{\thmname{#1}\thmnumber{ #2}\thmnote{ (#3)}}
 \theoremstyle{remark} 
\newtheorem{Remark}{Remark}

\usepackage{caption} 
\newcommand{\bound}{10^7}
\allowdisplaybreaks

\begin{document} 

\title[Numerical verification of Littlewood's bounds for $\vert L(1,\chi)\vert$]
{Numerical verification of Littlewood's bounds for $\vert L(1,\chi)\vert$}  
 \author{Alessandro Languasco} 

\begin{abstract}   
Let $L(s,\chi)$ be the Dirichlet $L$-function associated to a
non trivial primitive Dirichlet character $\chi$ defined $\bmod\ q$, where $q$ 
is an odd prime.
In this paper we introduce a fast method to compute $\vert L(1,\chi) \vert$  
using the values of Euler's $\Gamma$ function. 
We also introduce an alternative way of computing $\log \Gamma(x)$ and 
$\psi(x)= \Gamma^\prime/\Gamma(x)$,
$x\in(0,1)$. 
Using such algorithms we numerically verify the classical Littlewood  bounds
and the recent Lamzouri-Li-Soundararajan estimates on 
$\vert L(1,\chi) \vert$, where $\chi$  runs over the non trivial  primitive Dirichlet characters $\bmod\ q$,
 for every odd prime $q$ up to $\bound$.
The programs used and the results here described are collected 
at the following address \url{http://www.math.unipd.it/~languasc/Littlewood_ineq.html}.

\end{abstract} 
\maketitle
\makeatletter
\def\subsubsection{\@startsection{subsubsection}{3}%
  \z@{.3\linespacing\@plus.5\linespacing}{-.5em}%
  {\normalfont\bfseries}} 
\makeatother
\section{Introduction}
Let $q$ be an odd prime,  $\chi$ be a Dirichlet character $\bmod\ q$
and $L(s,\chi)$ be the associated Dirichlet $L$-function.
The goal of this paper is to  introduce a fast algorithm to compute 
the values of $\vert L(1,\chi) \vert$  
for every non trivial primitive Dirichlet character $\chi$ defined $\bmod\ q$
and, using such a new method, to 
numerically study a  generalisation of the 
classical bounds of Littlewood \cite{Littlewood1928}
for $\vert L(1,\chi_d) \vert$, where $\chi_d$ is a quadratic Dirichlet character.
Assuming the Riemann Hypotesis for $L(s,\chi_d)$ holds,
in 1928  Littlewood proved, for $d\ne m^2$, that
\begin{equation}
\label{Littlewood-bounds}
\Bigl(
\frac{12e^\gamma}{\pi^2}(1+\odi{1})\log \log \vert d\vert
\Bigr)^{-1}
<
L(1,\chi_d)
<
2e^\gamma (1+\odi{1}) \log \log \vert d\vert
\end{equation}
as $d$ tends to infinity, where $\gamma$ is the Euler-Mascheroni constant.
In 1973 Shanks \cite{Shanks1973} numerically studied 
the behaviour of the \emph{upper} and \emph{lower Littlewood indices} defined as
\begin{equation*}
\ULI (d,\chi_d):=
\frac{ L(1,\chi_d)}{2e^\gamma\log \log \vert d\vert}
\quad
\text{and}
\quad
\LLI(d,\chi_d) := L(1,\chi_d)\frac{12e^\gamma}{\pi^2}\log \log \vert d\vert
\end{equation*}
for several small discriminants $d$. Such computations were extended by 
Williams-Broere \cite{WilliamsB1976} in 1976 and by 
Jacobson-Ramachandran-Williams \cite{Jacobson2006} in 2006.

Recently  
Lamzouri-Li-Soundararajan \cite[Theorem 1.5]{LamzouriLS2015}
proved an effective form of Littlewood's inequalities:
assuming the Generalised Riemann Hypothesis holds, 
for every integer $q\ge 10^{10}$ and for every non trivial primitive character $\chi \bmod q$, 
they obtained that
\begin{equation}
\label{ULS-upper}
\vert L(1,\chi) \vert
\le 
2 e^\gamma \Bigl( \log \log q -\log 2 + \frac12 +\frac{1}{\log \log q} \Bigr)
\end{equation}
and
\begin{equation}
\label{LLS-lower}
\frac{1}{\vert L(1,\chi) \vert} 
\le 
\frac{12e^\gamma}{\pi^2}\Bigl( \log \log q - \log 2 + \frac12 +\frac{1}{\log \log q}  
+ \frac{14 \log \log q}{\log q}\Bigr).
\end{equation}

Using our  method we will compute the values of $\vert L(1,\chi) \vert$  
for every non trivial primitive Dirichlet character $\chi$ defined $\bmod\ q$, 
for every odd prime $q$ up to $\bound$.
This largely extends previous results.
Moreover, letting 
\begin{equation}
\label{Maxmin-def}
M_q:=\max_{\chi \neq \chi_0} 
\vert L(1,\chi) \vert, 
\quad 
m_q:=\min_{\chi \neq \chi_0} 
\vert L(1,\chi) \vert,
\end{equation}
\begin{equation}
\label{fg-def}
f(q):= \log \log q - \log 2 + 1/2 +1/\log \log q,
 \quad
g(q):= f(q) + 14 (\log \log q)/\log q,
\end{equation}
we obtain  the following
\begin{Theorem}
\label{Thm-max}
Let $3 \le q \le \bound$, $q$ be  a prime number and $M_q$ be defined in \eqref{Maxmin-def}.
 Then we have $0.604599\dotsc = M_3 \le M_q \le M_{4305479}= 6.399873\dotsc$. 
 Moreover, we also have
\[
0.325 \cdot 2 e^\gamma f(q)
<
M_q
< 
0.62 \cdot 2 e^\gamma f(q),
\]
where the lower bound holds just for $q\ge 79$, and
\[
0.4    
<
\max_{\chi\ne \chi_0 } \ULI(q,\chi)
<
0.66,  
\]
where the upper bound holds just for $q\ge 5$.
\end{Theorem}

 We also have an analogous result on $m_q$.
 \begin{Theorem}
 \label{Thm-min}
 Let $3 \le q \le \bound$, $q$ be  a prime number and $m_q$ be defined in \eqref{Maxmin-def}. Then
  we have $0.198814\dotsc = m_{991027} \le m_q \le m_{11} =0.618351\dotsc$. 
  Moreover, we also have
\[
 \frac{\pi^2}{12 e^\gamma} \frac{2.35}{g(q)}
<
m_q
<
  \frac{\pi^2}{12 e^\gamma} \frac{5}{g(q)},
\]
where the upper bound holds just for $q\ge 953$, and
\[
1.13
<
\min_{\chi\ne \chi_0 } \LLI(q,\chi)
<
2,
\]
where the lower bound holds just for $q\ge 373$.  
\end{Theorem}
 Theorems \ref{Thm-max}-\ref{Thm-min} are in agreement   with
 Littlewood's bounds in \eqref{Littlewood-bounds} and
the Lamzouri-Li-Soundararajan estimates
in \eqref{ULS-upper}-\eqref{LLS-lower}.  

The paper is organised as follows: in Section \ref{chi-Bernoulli-method} we will see how to compute 
$\vert L(1,\chi) \vert$ using the values 
of Euler's $\Gamma$ function and the Fast Fourier Transform algorithm;
we will also 
describe the actual computation we performed and how Theorems  \ref{Thm-max}-\ref{Thm-min} are obtained.
In Sections \ref{loggamma-computation}-\ref{psi-computation} we 
will see how to efficiently evaluate $\log\Gamma(x)$, and $\psi(x)=\Gamma^\prime/\Gamma(x)$, for $x\in (0,1)$ using 
precomputed values of the Riemann zeta-function at positive integers. After the bibliography 
we will also insert some tables and figures (the scatter plots were obtained using  GNUPLOT, v.5.2, patchlevel 8). 
 
 \medskip 
\textbf{Acknowledgements}. 
I wish to thank Luca Righi (University of Padova) for his help in developing the C language
implementation of the algorithms  described in Sections \ref{loggamma-computation}-\ref{psi-computation}.

\section{Computation of $\vert L(1,\chi)\vert$ and  proofs of Theorems \ref{Thm-max}-\ref{Thm-min}}
\label{chi-Bernoulli-method} 
Recall that $q$ is an odd prime and let $\chi$ be a primitive  non trivial Dirichlet character mod $q$.
The values of $\vert L(1,\chi) \vert$ can be computed in two different ways. 
Recalling  eq.~(3.1) of \cite{FordLM2014},   we have
$L(1,\chi) =-q^{-1}\sum_{a=1}^{q-1} \chi(a)\  \psi (a/q)$,
so that
\begin{equation}
\label{L-psi}
\vert L(1,\chi) \vert
=
\frac{1}{q}\
\Bigl \vert 
\sum_{a=1}^{q-1} \chi(a)\ 
\psi \bigl(\frac{a}{q}\bigr)
\Bigr\vert \ , 
\end{equation}
 where $\psi(x)= \Gamma^\prime/\Gamma(x)$ is the \emph{digamma} function
 and $\Gamma$ is Euler's function.
 As we will see later, for computational purposes it is in fact more 
 efficient to distinguish between the parity of the Dirichlet characters.
 If $\chi$ is an even character  we have, see, \emph{e.g.},
Cohen \cite[proof of Proposition 10.3.5]{Cohen2007}, that
 $L(1,\chi)
= 
2 \tau(\chi) q^{-1}
\sum_{a=1}^{q-1} \overline{\chi}(a)\log \bigl(\Gamma (a/q)\bigr)$, 
where the \emph{Gau\ss\ sum} 
$\tau(\chi):= \sum_{a=1}^q \chi(a)\, e(a/q)$, $e(x):=\exp(2\pi i x)$, verifies $\vert \tau(\chi) \vert = q^{1/2}$. Hence
\begin{equation}
\label{even} 
\vert 
L(1,\chi)
\vert 
= 
 \frac{2}{q^{1/2}}\
  \Bigl\vert 
\sum_{a=1}^{q-1} \overline{\chi}(a)\log\Bigl(\Gamma\bigl(\frac{a}{q}\bigr)\Bigr)\Bigr\vert
\quad 
(\chi\ \text{even}).
\end{equation}
Moreover, if $\chi$ is an odd character,  we have, see, \emph{e.g.},
Cohen \cite[Corollary 10.3.2]{Cohen2007}, that
  $L(1,\chi)
=
- w(\chi)  \pi  q^{-1/2} B_{1,\overline{\chi}}$, 
where $w(\chi) = \tau(\chi)/q^{1/2}$
and \(
B_{1,\chi}
:=
q^{-1} \sum_{a=1}^{q-1}  a \chi(a) 
\)
is the \emph{first $\chi$-Bernoulli number}.
 Hence  $\vert w(\chi)  \vert =1$ and
 \begin{equation}
\label{odd}
\vert 
L(1,\chi)
\vert 
=
  \frac{\pi}{q^{3/2}} \ \Bigl\vert \sum_{a=1}^{q-1}  a \overline{\chi}(a) \Bigr\vert
  \quad 
(\chi\ \text{odd}).
\end{equation}   
We will use the formulae \eqref{even}-\eqref{odd} because in half of the cases
we don't need  any special function, while in \eqref{L-psi}  
we need to evaluate the digamma function at $q-1$  points.
Moreover, in both the equations \eqref{even}-\eqref{odd} we can embed a \emph{decimation 
in frequency strategy} in the Fast Fourier Transform (FFT) algorithm  used to perform the sum over $a$,
see subsection \ref{fft-approach}.
Using the algorithm 
described in Section \ref{loggamma-computation}, see also Remark \ref{KaratsubaE-comparison}, 
the needed set of Gamma-function 
values can be computed with a precision of $n$ binary digits with a cost of 
$\Odi{qn}$ floating point products,  plus the cost of computing $(q-1)/2$ 
values of the logarithm function.
Hence, recalling also that the computational cost of the FFT algorithm 
of length $q$ is $\Odi{q\log q}$ floating point products, 
the total cost for computing $\vert L(1,\chi)\vert$ with a precision of $n$ binary digits
is then  $\Odi{q(n+\log q)}$ floating point products plus the cost of computing 
$(q-1)/2$ values of  the logarithm function. So far, this is
the fastest algorithm  to compute $\vert L(1,\chi)\vert$. 

We now proceed to describe our computational strategy.
Defining
\begin{align*} 
M^{\textrm{odd}}_q
&:=\max_{\chi\, \textrm{odd}}\
\vert L(1,\chi) \vert \ ,
\quad 
M^{\textrm{even}}_q
:=\max_{\substack{\chi \neq \chi_0\\ \chi\, \textrm{even}}} \
\vert L(1,\chi) \vert \ ,
\\ 
m^{\textrm{odd}}_q
&:=\min_{\chi\, \textrm{odd}}\
\vert L(1,\chi) \vert \ ,
\quad
m^{\textrm{even}}_q
:=\min_{\substack{\chi \neq \chi_0\\ \chi\, \textrm{even}}} \
\vert L(1,\chi) \vert \ ,
\end{align*}  
we will obtain $M_q,m_q$ as defined in \eqref{Maxmin-def}, using
\eqref{even}-\eqref{odd},
$M_q= \max(M^{\textrm{odd}}_q, M^{\textrm{even}}_q)$ and $m_q= \min(m^{\textrm{odd}}_q, m^{\textrm{even}}_q)$.

 \subsection{Computations trivially summing over $a$ (slower, more decimal digits available).}
In practice we first computed a few values of $M_q$ and $m_q$ using PARI/GP, v.~2.11.4, 
since it has  the ability to  generate
 the Dirichlet $L$-functions (and  many other $L$-functions).
 This can be done
with few instructions of the gp scripting language.
Such a computation has a linear cost in the number of calls
of the {\tt lfun} function of PARI/GP and it is, at least
on our Dell Optiplex desktop machine, slower than 
 using \eqref{even}-\eqref{odd}.
So we also implemented  such formulae in PARI/GP and
we were able to get the values of $M_q, m_q$ for every $q$ prime, $3\le q\le 1000$, 
with a precision of  $30$ decimal digits (see Tables \ref{table1} and \ref{table2})
in less than 17 seconds of computation time for each table. 
The machine we used was a Dell OptiPlex-3050, equipped with an Intel i5-7500 processor, 3.40GHz, 
16 GB of RAM and running Ubuntu 18.04.2.
   
  \subsection{Building the FFT approach}
\label{fft-approach} 
 As $q$ becomes large, the time spent in summing over $a$
dominates the overall computational cost.  So we implemented
the use of the FFT 
by using  the {\tt fftw} \cite{FFTW} library in our   C programs.   
We see now how to do so.

 In both   \eqref{even} and \eqref{odd} we remark that,
since $q$ is prime, it is enough to get $g$, a primitive root of $q$,
and $\chi_1$, the Dirichlet character mod $q$ given by 
 $\chi_1(g) = e^{2\pi i/(q-1)}$, to see that the set of the non-trivial characters
 mod $q$ is $\{\chi_1^j \colon j=1,2,\dotsc,q-2\}$.
 Hence, if, for every $k\in \{0,\dotsc,q-2\}$, we denote $g^k\equiv a_k\in\{1,\dotsc,q-1\}$,
 every summation in  \eqref{even}  and \eqref{odd}
 is  of the type $\sum_{k=0}^{q-2}  e^{-2\pi i j k /(q-1)} f(a_k/q)$, where $j\in\{1,\dotsc,q-2\}$
 is odd  and $f$ is a suitable function. 
As a consequence, such quantities are
 the Discrete Fourier Transform (DFT)  of the sequence $\{ f(a_k/q)\colon k=0,\dotsc,q-2\}$. 
This  idea was first formulated by Rader \cite{Rader1968} and it
 was used in \cite{FordLM2014,Languasco2019,LanguascoMSS2019,LanguascoR2020}
 to speed-up the computation of similar quantities via the use of Fast Fourier Transform dedicated
 software libraries. 
 
 In this case we can also use the \emph{decimation in frequency} strategy. Let $f$ be a
 function that assumes real values. Following the line in Section
 4.1 of \cite{Languasco2019}, letting $e(x):=\exp(2\pi i x)$, $m=(q-1)/2$,    
 for every $j=0,\dotsc, q-2$, $j=2t+\ell$, $\ell\in\{0,1\}$ and $t\in \Z$, we have that  
\[
 \sum_{k=0}^{q-2}   e\Bigl(\frac{- j k}{q-1}\Bigr)  f \Bigl(\frac{a_k}{q}\Bigr)
 =
 \sum_{k=0}^{m-1}  
 e\Bigl(\frac{- t k}{m}\Bigr)  
    e\Bigl(\frac{- \ell k}{q-1}\Bigr)    
 \Bigl(
 f\Bigl(\frac{a_k}{q}\Bigr) 
 +
 (-1)^{\ell} 
 f \Bigl(\frac{a_{k+m}}{q}\Bigr)
 \Bigr),
\]
 where $t=0,\dotsc, m-1$. 
 Letting
 \begin {equation} 
\label{bk-ck-def} 
b_k :=
  f\Bigl(\frac{a_k}{q}\Bigr) +  f \Bigl(\frac{a_{k+m}}{q}\Bigr)   
\quad
\textrm{and}
\quad
c_k :=  
   e\Bigl(-\frac{k}{q-1}\Bigr)   
 \Bigl(  f\Bigl(\frac{a_k}{q}\Bigr) -  f \Bigl(\frac{a_{k+m}}{q}\Bigr)  \Bigr),
\end{equation}
we can rewrite the previous formula  (recall that $j=2t+\ell$, $\ell\in\{0,1\}$ and $t=0,\dotsc, m-1$) as
 \begin {equation} 
\label{DIF} 
   \sum_{k=0}^{q-2}   e\Bigl(\frac{- j k}{q-1}\Bigr)  f \Bigl(\frac{a_k}{q}\Bigr)
 =
 \begin{cases}
 \sum\limits_{k=0}^{m-1}     e\bigl(-\frac{t k}{m}\bigr) b_k 
 & \textrm{if} \ \ell =0\\
  \sum\limits_{k=0}^{m-1}     e\bigl(-\frac{t k}{m}\bigr)  c_k  
 & \textrm{if} \ \ell =1.\\
 \end{cases}
\end{equation}
Since we just need the sum  over the odd Dirichlet characters for $f(x)=x$ 
and over the even Dirichlet characters for $f(x)=\log \Gamma(x)$,
in this way we can evaluate  an FFT of length $(q-1)/2$, instead of $q-1$,
applied on a suitably modified sequence according to \eqref{bk-ck-def}-\eqref{DIF}.
Clearly this represents a gain in both speed and memory usage
in running the actual computer program. 

In the case $f(x)= \log \Gamma(x)$ we can simplify the form of 
$b_k = \log \Gamma(a_k/q) + \log \Gamma(a_{k+m}/q) $, where $m=(q-1)/2$ and $k=0,\dotsc, m-1$,
in the following way. Recalling $\langle g \rangle = \Z^*_q$, $a_k \equiv g^k \bmod q$ and   $g^m \equiv q-1 \bmod{q}$,
 we can write that  $ a_{k+m}  \equiv  q-a_{k} \bmod{q}$
 and hence
\(
\log \Gamma(a_{k+m}/q)  = \log \Gamma ((q-a_{k})/q)
=
\log \Gamma (1- a_{k}/q).
\)
Using the well-known \emph{reflection formula} $\Gamma(x) \Gamma(1-x)  = \pi / \sin(\pi x)$,  
we obtain 
\begin{align}
\label{DIF-sin}
\log \Gamma\Bigl(\frac{a_{k}}{q}\Bigr)
+
\log \Gamma\Bigl(\frac{a_{k+m}}{q}\Bigr)
=
\log \Gamma\Bigl(\frac{a_{k}}{q}\Bigr)
+
\log \Gamma\Bigl(1-\frac{a_{k}}{q}\Bigr)
&=
 \log \pi - \log \Bigl(\sin\bigl( \frac{\pi a_k}{q}\bigr)\Bigr) ,
\end{align}
for every $k=0,\dotsc, m-1$. Inserting the last relation in the definition of $b_k$  in \eqref{bk-ck-def}
and remarking that, by orthogonality, the constant term $\log \pi$ is negligible,
we can replace in the actual computation  the  Gamma function  with the $\log(\sin (\cdot))$ one.
Since in our application we will have $a/q\in(0,1)$,  
we also developed our own alternative implementation of $\log \Gamma(x)$, $x\in (0,1)$, 
see Section \ref{loggamma-computation}. 

In the case $f(x)=x$, it is easier to obtain a simplified  form of $c_k$  as defined in \eqref{bk-ck-def}. 
Using again $\langle g \rangle = \Z^*_q$, 
$a_k \equiv g^k \bmod q$ and   $g^m \equiv q-1 \bmod{q}$,
 we can write that $ a_{k+m}  \equiv  q-a_{k} \bmod{q}$;  hence
\(
 a_k  -    a_{k+m}  
 = 
  a_k -(q-a_{k})  
=
2a_k  -q,
\)
so that in this case,  for every $k=0,\dotsc m-1$, $m=(q-1)/2$, we obtain 
\[
c_k=   e\Bigl(-\frac{k}{q-1}\Bigr)\Bigl(2\frac{a_k}{q} -1\Bigr).
\]

\subsection{Computations  summing over $a$ via FFT (much faster, less decimal digits available).}
\label{practical-comp}
\begin{figure}[ht] 
\includegraphics[scale=0.48,angle=0]{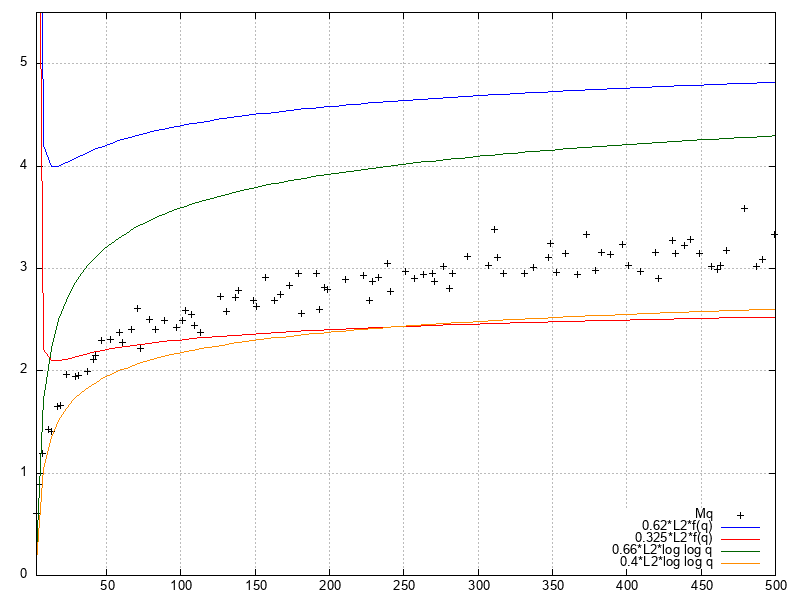}  
\caption{\small{Graphical representation of Theorem \ref{Thm-max}.
The values of $M_q$, $q$ prime, $3\le q\le 500$.  
The blue  line  represents $ 0.62 L_2  f(q)$;  the red one  
$0.325 L_2  f(q)$, where   $f(q)$ is defined in \eqref{fg-def}
and $L_2=  2 e^\gamma$.
The green line represents $0.66L_2\log \log q$; the orange one $0.4L_2\log \log q$.
}
 }
\label{fig-LLS-1-small} 
 \end{figure}

\begin{figure}[ht] 
\includegraphics[scale=0.5,angle=0]{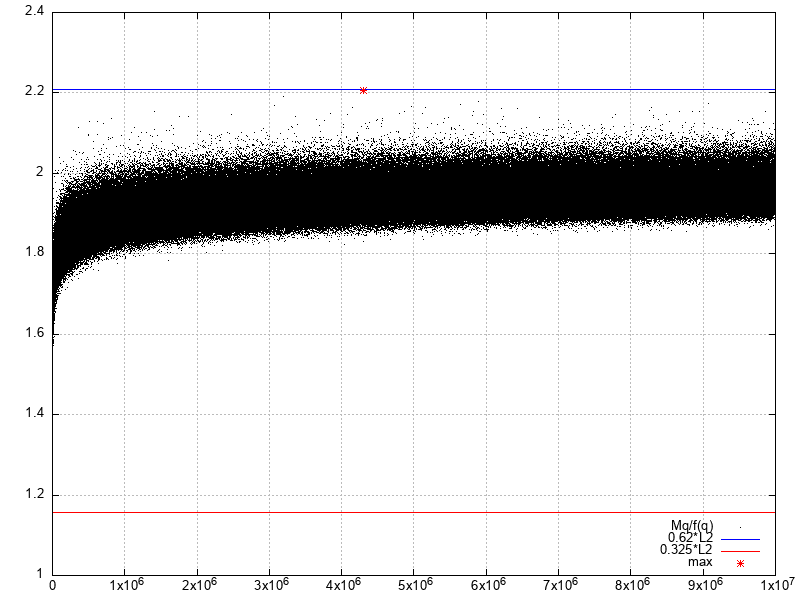}  
\caption{\small{The values of $M_q^\prime:=M_q/f(q)$, $q$ prime, $500\le q\le \bound$, 
where $f(q)$ is defined in \eqref{fg-def}. 
The minimal value for $M_q^\prime$ is $0.057396\dotsc$  attained at $q=3$ and the maximal one is  
$2.206927\dotsc$ attained at $q=4305479$ (the ``second'' maximal value is  $2.192260\dotsc$
attained at $q=3190151$).
The blue  line  represents $ 0.62 L_2$;  the red one  
$0.325 L_2$,  where  $L_2=  2 e^\gamma$.
For $3\le q <500$, see the next plot.
 }
 }
\label{fig-LLS-2} 
 \end{figure}

\begin{figure}[ht] 
\includegraphics[scale=0.5,angle=0]{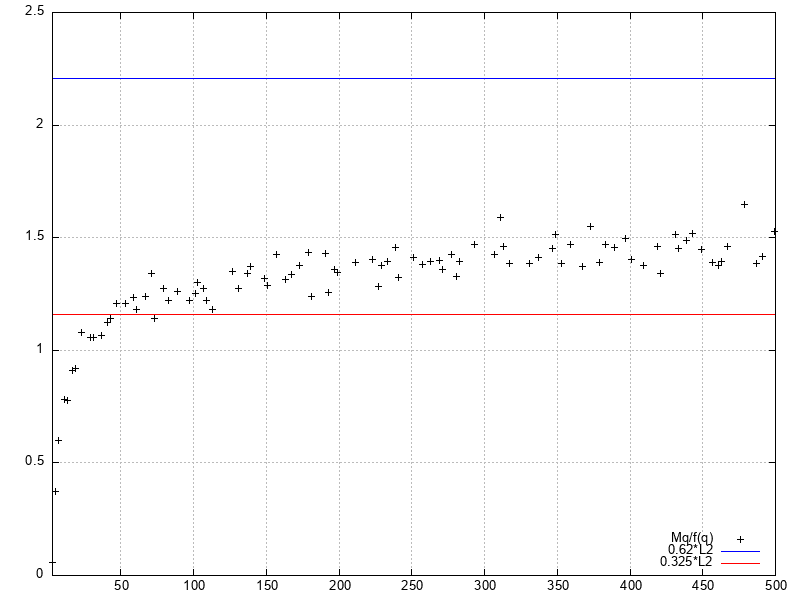}  
\caption{\small{The values of $M_q^\prime:=M_q/f(q)$, $q$ prime, $3\le q\le 500$, where $f(q)$ 
is defined in \eqref{fg-def}. 
The blue  line  represents $ 0.62 L_2$;  the red one  
$0.325 L_2$,  where  $L_2=  2 e^\gamma$.
 }
 }
\label{fig-LLS-2-small}  
 \end{figure}

\begin{figure}[ht] 
\includegraphics[scale=0.5,angle=0]{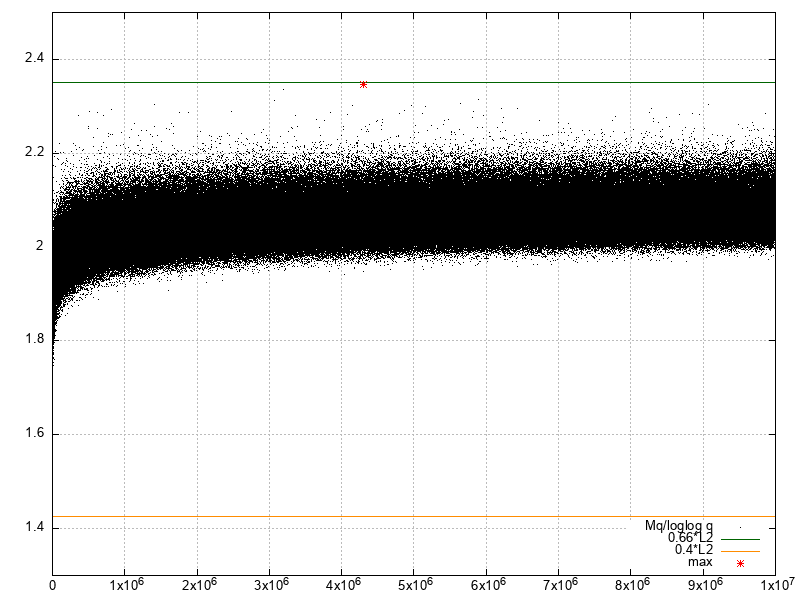}  
\caption{\small{The values of $M_q^{\prime\prime}:=M_q/\log \log q$, $q$ prime, 
$500\le q\le \bound$.  
The minimal value for $M_q^{\prime\prime}$ is $1.492809\dotsc$  attained at $q=13$ and the maximal one is   
$6.428641 \dotsc$ attained at $q= 3$ (the ``second'' maximal value is  $2.347506\dotsc$
attained at $q=4305479$, marked in red in this plot).
The green  line  represents $ 0.66 L_2$;  the orange one  
$0.4 L_2$,  where  $L_2=  2 e^\gamma$.
For $3\le q <500$, see the next plot.
 }
 }
\label{fig-ULI} 
\end{figure}

\begin{figure}[ht] 
\includegraphics[scale=0.5,angle=0]{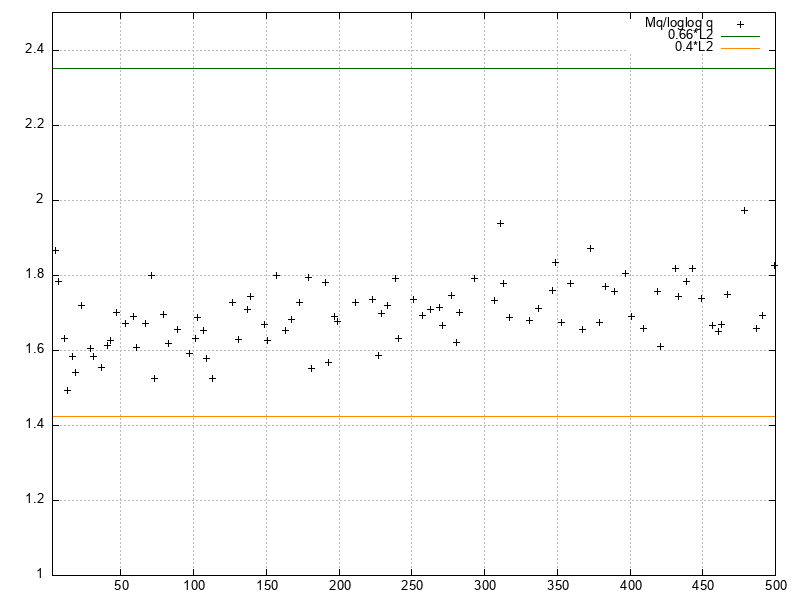}  
\caption{\small{The values of $M_q^{\prime\prime}:=M_q/\log \log q$, $q$ prime, $5\le q\le 500$. 
The green  line  represents $ 0.66 L_2$;  the orange one  
$0.4 L_2$,  where  $L_2=  2 e^\gamma$.
$M_3^{\prime\prime}=6.428641\dots$ is not included in this plot.
 }
 }
\label{fig-ULI-small} 
 \end{figure}

\begin{figure}[ht] 
\includegraphics[scale=0.48,angle=0]{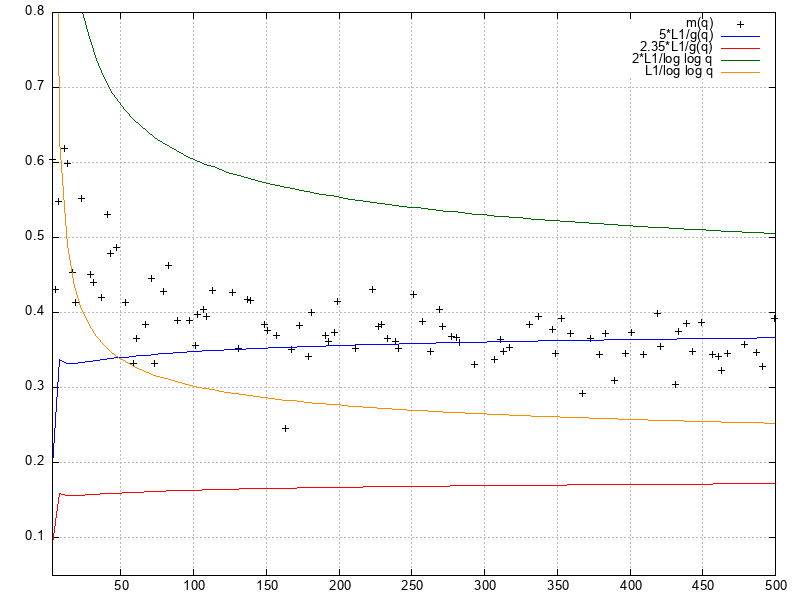}  
\caption{\small{Graphical representation of Theorem \ref{Thm-min}.
The values of $m_q$, $q$ prime, $500\le q\le \bound$. 
The blue  line  represents $ 5 L_1 / g(q)$;  the red one  
$2.35 L_1 / g(q)$, where $g(q)$ is defined in \eqref{fg-def}
and  $L_1=   \frac{\pi^2}{12 e^\gamma}$.
The green line represents $2L_1/\log\log q$; the orange one $1.13 L_1/\log\log q$.
}
 }
\label{fig-LLS-3-small} 
 \end{figure}

\begin{figure}[ht] 
\includegraphics[scale=0.5,angle=0]{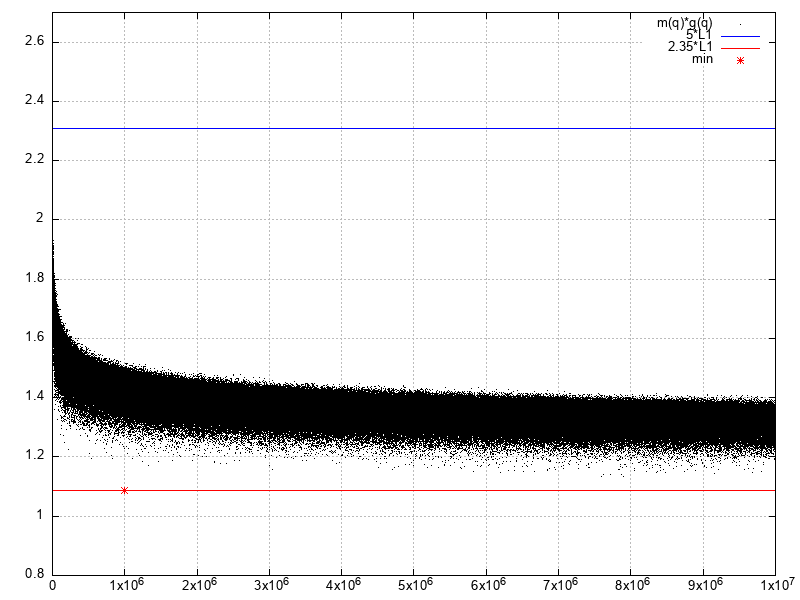}  
\caption{\small{The values of $m_q^\prime:=m_q g(q)$, $q$ prime, $500\le q\le \bound$, where $g(q)$ is defined in \eqref{fg-def}.  
The minimal value for $m_q^\prime$ is $1.088477\dotsc$ attained at $q=991027$
(the ``second'' minimal value is  $1.134017\dotsc$
attained at $q=7598287$)
and the maximal one is   $7.093329\dotsc$ attained at $q=3$.
The blue  line  represents $ 5 L_1 $;  the red one  
$2.35 L_1 $, where  $L_1=   \frac{\pi^2}{12 e^\gamma}$.
 For $3\le q <500$, see the next plot.
 }
 }
\label{fig-LLS-4} 
 \end{figure}

\begin{figure}[ht] 
\includegraphics[scale=0.5,angle=0]{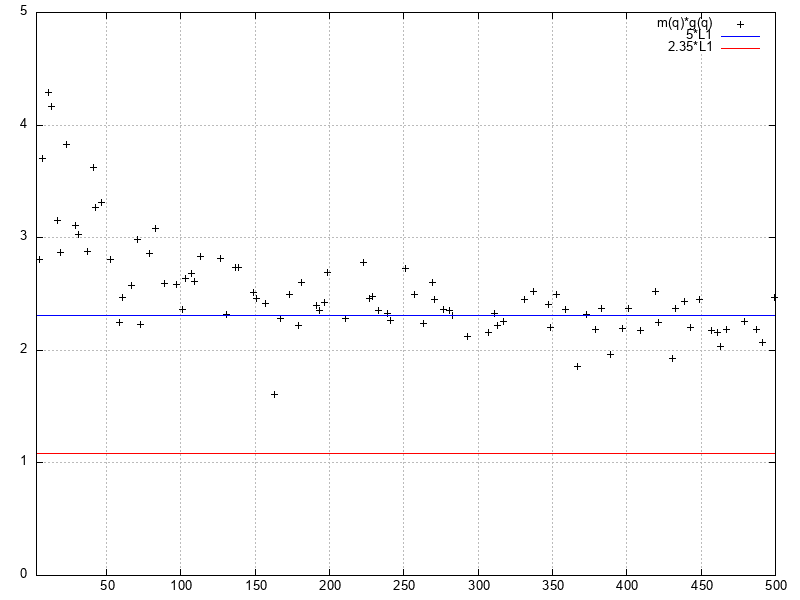}  
\caption{\small{The values of $m_q^\prime:=m_q g(q)$, $q$ prime, $3\le q\le 500$, 
where $g(q)$ is defined in \eqref{fg-def}.  
The blue  line  represents $ 5 L_1 $;  the red one  
$2.35 L_1 $, where  $L_1=   \frac{\pi^2}{12 e^\gamma}$. 
$m_3^{\prime}=7.093329\dotsc$ is not included in this plot.
 }
 }
\label{fig-LLS-4-small} 
 \end{figure}

\begin{figure}[ht] 
\includegraphics[scale=0.5,angle=0]{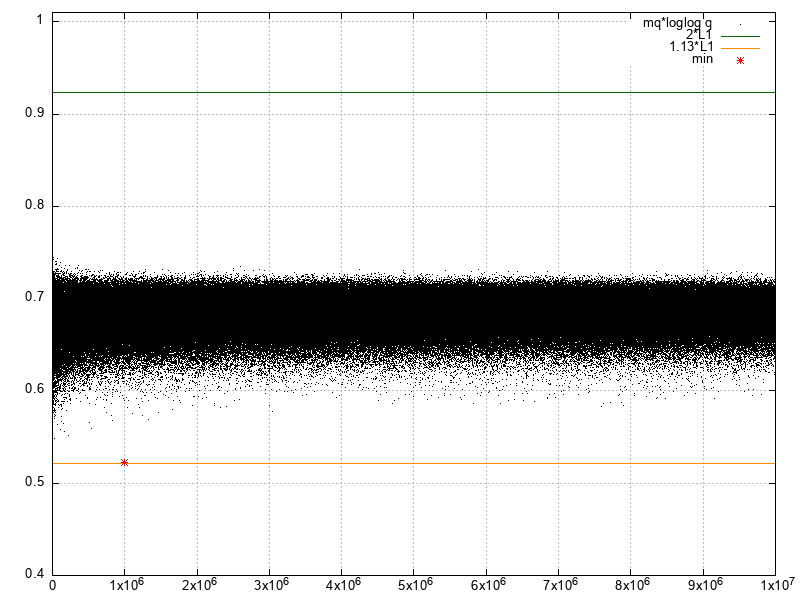}  
\caption{\small{The values of $m_q^{\prime\prime}:=m_q\log \log q$, 
$q$ prime, $500\le q\le \bound$.  
The minimal value for $m_q^{\prime\prime}$ is $0.056861\dotsc$ attained at $q=3$
(the ``second'' minimal value is  $0.2048251\dotsc$
attained at $q=5$)
and the maximal one is   $0.7445135\dotsc$ attained at $q=19001$. 
The green  line  represents $ 2L_1$;  the orange one  
$1.13 L_1$,  where   $L_1=   \frac{\pi^2}{12 e^\gamma}$. 
The red point represents $m_{991027} ^{\prime\prime}= 0.521914\dots$
For $3\le q <500$, see the next plot.
 }
 } 
\label{fig-LLI} 
\end{figure}

\begin{figure}[ht] 
\includegraphics[scale=0.5,angle=0]{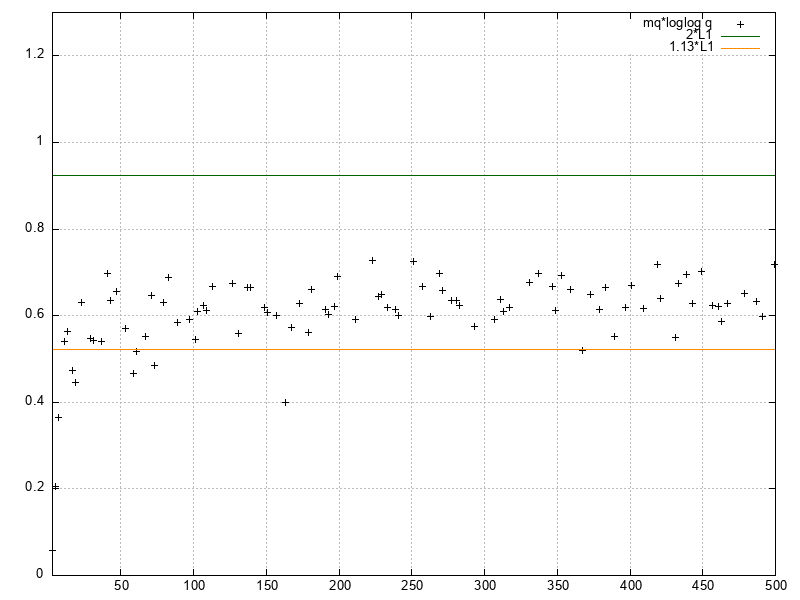}  
\caption{\small{The values of $m_q^{\prime\prime}:=m_q\log \log q$, $q$ prime. 
The green  line  represents $ 2L_1$;  the orange one  
$1.13 L_1$,  where   $L_1=   \frac{\pi^2}{12 e^\gamma}$. 
 }
 }
\label{fig-LLI-small} 
 \end{figure}

Using the setting explained in the previous subsection,  we were able to compute, 
using the \emph{long double precision} (80 bits)
of the C programming language, the values of $M_q$ and $m_q$ 
for every prime $3\le q \le  \bound$ and we provide here 
the scatter plots of such values and of their normalisations, see Figures \ref{fig-LLS-1}-\ref{fig-LLI-small}.
The data were obtained  in about 
$57$ days  of computation time on 
the Dell OptiPlex machine mentioned before. 
  
The actual FFTs were performed using the FFTW \cite{FFTW} software library.
The PARI/GP scripts and the C programs used and the computational results obtained
are available at the following web address:
\url{http://www.math.unipd.it/~languasc/Littlewood_ineq.html}. 

\subsection{Proof of Theorems \ref{Thm-max}-\ref{Thm-min}}  
Theorems \ref{Thm-max}-\ref{Thm-min} follow by analysing,
using  suitable programs  written in \texttt{python},  the
data computed in subsection \ref{practical-comp} and collected 
in two comma-separated values (csv) files. 
We obtain that the inequalities in the statements of Theorems 1-2 hold
and that the \newcommand{\captionmax}{minimal value for $M_q$ is $0.604599\dotsc$  attained at $q=3$ and the maximal one
is  $6.399873\dotsc$ attained at $q=4305479$. 
}
\captionmax
\newcommand{\captionmin}{The minimal value for $m_q$ is $0.198814\dotsc$ attained at $q=991027$ 
 and the maximal one is   $0.618351\dotsc$ attained at $q=11$.}
 \captionmin  

The output of such \texttt{python} programs are available at the web page:
\url{http://www.math.unipd.it/~languasc/Littlewood_ineq.html}.
A few plots representing  Theorems \ref{Thm-max}-\ref{Thm-min} 
are given in Figures \ref{fig-LLS-1}-\ref{fig-LLS-1-small} and \ref{fig-LLS-3}-\ref{fig-LLS-3-small}.

\section{An alternative algorithm to compute $\log \Gamma(x)$, $x\in(0,1)$}
\label{loggamma-computation}

We describe here an alternative way of computing 
$\log \Gamma(x)$, $x\in(0,1)$,
which is based on the well-known Euler formula (see, \emph{e.g.}, Lagarias 
\cite[section 3]{Lagarias2013}):
\begin{equation}
\label{Euler-loggamma}
\log\Gamma(x) = \gamma(1-x)  + \sum_{k=2}^{+\infty} \frac{\zeta(k)}{k}(1-x)^k,
\end{equation}
where $\zeta(s)$ is the Riemann zeta-function.

We follow the argument used in Languasco-Righi \cite{LanguascoR2020} to study
the Ramanujan-Deninger Gamma function $\Gamma_1(x)$.
We immediately remark that the series in \eqref{Euler-loggamma}
absolutely converges for $x\in(0,2)$; this fact and the well-known relation 
\begin{equation}
\label{difference-gamma}
\log \Gamma (1+x) =\log \Gamma(x)  + \log x, \quad x>0,
\end{equation}
 let us 
obtain $\log \Gamma(x)$,  $x\in (0,1)$, in two different ways.
Recalling $\log \Gamma(1)=0$ and $\log \Gamma(1/2)=(\log \pi)  /2$,
we also remark that, letting $n\in \N$, $n\ge 2$, 
for every $x\in(0,2)$
there exists $r=r_\Gamma(x,n)  \ge 2$ 
such that
\begin{equation}
\label{tail-euler-series}
\Bigl\vert \sum_{k=r+1}^{+\infty} \frac{\zeta(k)}{k}(1-x)^k \Bigr\vert
<
\frac{\zeta(3)}{3}
\sum_{k=r+1}^{+\infty} 
\vert 1-x \vert^{k} 
<
0.41
\frac{\vert 1-x \vert^{r+1}}{1-\vert 1-x\vert}
<2^{-n-1}.
\end{equation}
A straightforward computation reveals that we can choose
\begin{equation*}
r_\Gamma(x,n) = 
\Bigl\lceil \frac{(n+1) \log 2 + \vert \log (1-\vert 1-x\vert)\vert}{ \vert \log  \vert 1-x \vert \vert}\Bigr\rceil -1,
\end{equation*}
where we denoted  as $\lceil y \rceil$ 
the least integer greater than or equal to 
$y\in \R$.

\subsection{The shifting trick for $\log \Gamma(x)$, $x\in (0,1)$}
\label{shifting-trick}
Clearly $r_\Gamma(x,n)$ becomes larger as $\vert 1-x\vert$ increases. So
when $x$ is close to zero we will   evaluate $\log \Gamma$ at $1+x$
via \eqref{difference-gamma}.
In the following we will refer to this idea as the \emph{shifting trick}.
This way we will always use the best convergence interval, $x\in(1/2,3/2)$, we have for
the series in \eqref{Euler-loggamma}; we also remark that 
$r_\Gamma(x,n) \le r_\Gamma(1/2,n) =r_\Gamma(3/2,n) =  n+1$ for every $x\in(1/2,3/2)$.  
Summarising,  using  \eqref{Euler-loggamma} and \eqref{tail-euler-series},
 for $x\in(1/2,1)$ we have that there exists $\theta=\theta(x) \in (-1/2,1/2)$ such that
 \begin{equation}
 \label{Gamma-x>1/2}
\log \Gamma(x)
=
 \gamma (1-x)  + \sum_{k=2}^{+\infty} \frac{\zeta(k)}{k}(1-x)^k
=
 \gamma (1-x)  + \sum_{k=2}^{r_\Gamma(x,n)} \frac{\zeta(k)}{k}(1-x)^k
 +\vert \theta \vert 2^{-n}.
  \end{equation}
  We also  remark that for  $x\in(1/2,1)$, we have 
\begin{equation*}
r_\Gamma(x,n) = \Bigl\lceil \frac{(n+1) \log 2 + \vert \log x\vert}{ \vert \log  ( 1-x ) \vert} \Bigr\rceil -1
\le n+1.
\end{equation*}
Moreover, using  \eqref{Euler-loggamma}-\eqref{tail-euler-series}, 
for $x\in(0,1/2)$ we have  that there exists $\eta=\eta(x) \in (-1/2,1/2)$ such that
 \begin{equation}
 \label{Gamma-x<1/2}
\log \Gamma(x)
= 
- \log x  -\gamma x  + \sum_{k=2}^{+\infty} \frac{(-1)^k \zeta(k)}{k} x^k
=
- \log x  -\gamma x  + \sum_{k=2}^{r^\prime_\Gamma(x,n)} \frac{(-1)^k \zeta(k)}{k} x^k
+\vert \eta \vert 2^{-n},  
  \end{equation}
  where 
\begin{equation*}
r^\prime_\Gamma(x,n):=r_\Gamma(1+x,n) 
  = \Bigl\lceil \frac{(n+1) \log 2 + \vert \log (1-  x )\vert}{ \vert \log    x   \vert} \Bigr\rceil-1
  \le n+1.
\end{equation*} 
Since the needed $\zeta$-values  can be precomputed and stored with the desired precision
(using, for example, PARI/GP),
 the formulae in \eqref{Gamma-x>1/2}-\eqref{Gamma-x<1/2} allow us to compute
$\log \Gamma(x)$, $x\in (0,1)$, with a precision of $n$ binary digits 
using at most $n+1$ summands; moreover, they also reveal that  such a task is, 
from a computational point of view,
 essentially as difficult
as computing   $\log x$ when $x$ is close to $0$.

\begin{Remark}[Computational cost]
\label{KaratsubaE-comparison}
The estimates $r_\Gamma(x,n),r^\prime_\Gamma(x,n) \le  n+1$ for
every $x\in (1/2,1)$ and, respectively, $x\in (0,1/2)$, imply
that $\log \Gamma(x)$, $x\in (0,1)$ can be obtained with a $n$-bit
precision using at most $n+1$ summands.
The summation is performed combining 
the ``pairwise summation'' \cite{Higham1993}  algorithm
with Kahan's \cite{Kahan1965} method  (the minimal block
for the pairwise summation algorithm is summed using 
Kahan's method)
to have a good compromise between precision, computational
cost and execution speed.
Hence the cost of computing $\log \Gamma(x)$, $x\in (1/2,1)$
is  $\Odi{n}$ floating point products and $\Odi{n}$ floating point summations
with a precision of $n$ binary digits; for $x\in (0,1/2)$ we have the same plus the cost of computing $\log x$.

In the particular case in which $x=a/q$ and $a$ runs over $1,\dotsc, q-1$,
the total cost to obtain the values  $\{\log \Gamma(a/q)\colon  a=1, \dotsc, q-1\}$, each one 
with a precision of $n$ binary digits, is  then $\Odi{qn}$  floating point products, 
plus the cost of computing $(q-1)/2$ 
values of the logarithm function.
\end{Remark}

\begin{Remark}[Computation in the whole real axis]
It is clear that using \eqref{difference-gamma} and \eqref{Gamma-x>1/2}-\eqref{Gamma-x<1/2}
we can compute $\log \Gamma(x)$ for every $x>0$ as follows.
For every $x>0$, we denote
as $\lfloor x \rfloor$ the integral part of $x$ and as $\{x\} = x - \lfloor x \rfloor$
the fractional part of $x$. Hence we obtain:
\begin{enumerate}[i)]
\item
 $\log \Gamma(1)=\log \Gamma(2)=0$ and   $\log\Gamma(m)  = \sum_{k=2}^{m-1} \log k$
for every $m\in \N$, $m\ge 3$;
\item
for  $x>1$, $x\not \in \N$, we  compute $\log\Gamma(x)$ as $\log\Gamma(x) 
= \log\Gamma(\{x\}) + \sum_{k=0}^{\lfloor x \rfloor - 1} \log (\{x\}+k)$;
\item
$\log\Gamma(1/2) =   (\log \pi)/2$;
\item
for $x\in (0,1/2)$,  we  compute $\log\Gamma(x)$ as in \eqref{Gamma-x<1/2};
\item
for  $x\in (1/2,1)$,   we  compute $\log\Gamma(x)$ as in \eqref{Gamma-x>1/2}.
\end{enumerate}
Even if we are mainly interested in working with $x\in(0,1)$ we recall that for $x$ large it 
might be more convenient to implement Stirling's formula for $\log \Gamma(x)$.
\end{Remark}

\begin{Remark} [Enlarging the convergence radius] We remark, even if it is not 
useful in our application, that the size of the convergence interval  in 
\eqref{Euler-loggamma} can be doubled by isolating
the Taylor series at $1$ of $\log x- (x-1)$ in \eqref{Euler-loggamma}  thus getting
 
\begin{equation}
\label{enlarged-radius}
\log\Gamma(x) = - \log x+ ( \gamma-1)(1-x)  + \sum_{k=2}^{+\infty} \frac{\zeta(k)-1}{k}(1-x)^k.
\end{equation}
Using the well-known estimate
$\vert \zeta(k)-1 \vert<2^{1-k}$ for every $k\in\N, k\ge 3$,
it is easy to prove that  the  series in \eqref{enlarged-radius} converges for every $x\in(-1,3)$.
\end{Remark}

\begin{Remark}[Computation in the complex plane]
\begin{enumerate}[1)]
\item
It seems that the argument leading to \eqref{Gamma-x>1/2}-\eqref{Gamma-x<1/2}
is not usually implemented in the most used software libraries or Computer Aided Systems (CAS) for Mathematics
probably because the shifting trick used before can be directly generalised to complex variables
only in a thin horizontal strip around the positive part of the real axis, see the next point
of this remark. 
In fact, many software libraries and CAS usually implement the 
computation of $\log \Gamma (z)$, $z\in \C$, $z\ne -n$, $n\in \N$, using the Lanczos 
approximation thus following the setting of  Press \emph{et al.} 
\cite{Press1992}.
\item

A possible complex strip can be built 
combining  \eqref{difference-gamma}, which in fact holds for any argument 
$z\in \C$, $z\ne -n$, $n\in \N$, together with
 the complex power series contained in the following formula

\begin{equation}
\label{loggamma-complex}
\log\Gamma(z) 
= 
 \gamma(1-z)  + \sum_{k=2}^{+\infty} \frac{\zeta(k)}{k}(1-z)^k
\end{equation}
which generalises \eqref{Euler-loggamma} to the region   $\vert z-1\vert<1$.
We can start from the rectangle 
$1/2\le \Re(z)< 3/2$, $\vert \Im(z)\vert \le 1/4$,
since for every $z$ in this region less than 
$1.2\cdot(n+1)+4$ terms
are sufficient to have a precision of $n$ binary digits in computing a truncation of the series
in \eqref{loggamma-complex}. 
We also remark here that such a strip can be  vertically
enlarged using Gau\ss' multiplication
theorem in the following form
\[
\log \Gamma(mz) =
\frac{1-m}{2}\log(2\pi) 
+
\Bigl(mz-\frac12\Bigr)\log m 
+ 
\sum_{j=0}^{m-1}\log \Gamma \Bigl(z+\frac{j}{m}\Bigr),
\]
where $m\in \N$, $m\ge 1$.
\end{enumerate}
\end{Remark}

\subsection{Reflection formulae} 
We now remark that
using    \eqref{Gamma-x>1/2}-\eqref{Gamma-x<1/2} to compute $\log\Gamma(x) +\log\Gamma(1-x)$,
$x\in (0,1)$,   the odd summands of the series will vanish and something similar happens in computing
$\log\Gamma(x) -\log\Gamma(1-x)$.
We summarise the situation in the following
\begin{Proposition}
\label{DIF-formulae-gamma}
Let $x\in (0,1)$, $x\ne 1/2$,
 $n\in \N$, $n\ge 2$, $r_1(x,n) =
 \lceil\frac{(n+1) \log 2 + \vert \log (1-x) \vert}{ \vert \log  x \vert } -1\rceil/2$
and $r_2(x,n) = \lceil\frac{(n+1) \log 2 +\vert \log x \vert}{\vert \log  (1-x) \vert}-1\rceil/2$.
Using  \eqref{Gamma-x>1/2} and \eqref{Gamma-x<1/2},
we have that there exists $\theta=\theta(x) \in (-1/2,1/2)$  such that
for $0<x  < 1/2$ we have
\begin{align} 
\log\Gamma(x)&+\log\Gamma(1-x) 
 \label{Gamma-DIF-x<1/2}
= 
-  \log x + 
\sum_{\ell=1}^{r_1} \frac{\zeta(2\ell)}{\ell} x^{2\ell}
+ \vert \theta \vert 2^{-n}, 
\\
\log\Gamma(x)&-\log\Gamma(1-x) 
 \label{Gamma-DIF-odd-x<1/2}
= 
-  \log x  - 2\gamma x 
-2 
\sum_{\ell=1}^{r_1} \frac{\zeta(2\ell+1)}{2\ell+1} x^{2\ell+1}
+ \vert \theta \vert 2^{-n}, 
\end{align}
and  for $1/2<x<1$ we have
\begin{align}
\log\Gamma(x)&+\log\Gamma(1-x) 
\label{Gamma-DIF-x>1/2}
=
-  \log (1-x) + 
\sum_{\ell=1}^{r_2} \frac{\zeta(2\ell)}{\ell} (1-x)^{2\ell}
+ \vert \theta \vert 2^{-n},
 \\ 
\log\Gamma(x)&-\log\Gamma(1-x) 
\label{Gamma-DIF-odd-x>1/2}
=
 \log (1-x) +  2\gamma(1-x)
+2 \sum_{\ell=1}^{r_2} \frac{\zeta(2\ell+1)}{2\ell+1} (1-x)^{2\ell+1}
+ \vert \theta \vert 2^{-n}.
 \end{align}
\end{Proposition}
 \begin{Proof}
Assume that $0<x < 1/2$; we compute $\log\Gamma(x)$ with the infinite 
series in  \eqref{Gamma-x<1/2} 
and $\log\Gamma(1-x)$ with the infinite series in  \eqref{Gamma-x>1/2}. 
Since they absolutely converge, their sum is obtained with the series
having as summands the sum of their coefficients. Arguing as in 
\eqref{tail-euler-series} and remarking that  
$r_1(x,n)= r_\Gamma(1-x,n) /2=  r^\prime_\Gamma(x,n)/2$,
we immediately have  that \eqref{Gamma-DIF-x<1/2} holds since
  the odd summands vanish.
Assume that $1/2<x<1$; in this case we  compute $\log\Gamma(x)$ with the 
infinite series in \eqref{Gamma-x>1/2} 
and $\log\Gamma(1-x)$ with the infinite series in  \eqref{Gamma-x<1/2}. 
Since they absolutely converge, their sum is obtained with the series
having as summands the sum of their coefficients. Arguing as in 
\eqref{tail-euler-series} and remarking that  
$r_2(x,n)= r_\Gamma(x,n)/2 =  r^\prime_\Gamma(1-x,n)/2$, we
immediately have  that \eqref{Gamma-DIF-x>1/2} holds since
  the odd summands vanish. 
 The  derivation of \eqref{Gamma-DIF-odd-x<1/2} and \eqref{Gamma-DIF-odd-x>1/2}
 is similar.
  This completes the proof.
\end{Proof}
 
It is worth mentioning that the right hand side in \eqref{Gamma-DIF-x>1/2} can
be obtained from the one in \eqref{Gamma-DIF-x<1/2} formally replacing $x$ with $1-x$;
and that, also changing of sign, the same holds for  \eqref{Gamma-DIF-odd-x<1/2}    
and \eqref{Gamma-DIF-odd-x>1/2}.
Using $\Gamma(x) \Gamma(1-x)  = \pi / \sin(\pi x)$,   Proposition \ref{DIF-formulae-gamma}
 also immediately gives a way of writing $\sin(\pi x)$
in term of logs and values of the Riemann zeta-function at positive even integers,
see Remark \ref{sin-zeta} below.  
Comparing with \eqref{DIF-sin},
the use of Proposition \ref{DIF-formulae-gamma} in our application 
is particularly efficient for the following reasons:
\begin{enumerate}[$\bullet$] 

 \item the cancellation of the odd terms we have in \eqref{Gamma-DIF-x<1/2} and \eqref{Gamma-DIF-x>1/2}
 leads to  gain  a factor  of $2$  in the computational cost
  since we just need to use half of the summands (the ones with even indices);
  a similar remark applies to \eqref{Gamma-DIF-odd-x<1/2} and \eqref{Gamma-DIF-odd-x>1/2} too;
  
  \item in \eqref{Gamma-DIF-x<1/2} and \eqref{Gamma-DIF-x>1/2} the  number of summands
  is $\le (n+1)/2$; hence to have a precision of $n$ bits we just need less than $(n+1)/2$
  summands (assuming the logarithm function can be evaluated with the same precision);
    a similar remark applies to \eqref{Gamma-DIF-odd-x<1/2} and \eqref{Gamma-DIF-odd-x>1/2} too;
  
 \item in \eqref{Gamma-DIF-x<1/2} and \eqref{Gamma-DIF-x>1/2} just the values of
the Riemann zeta-function at positive even integers are required and for them
we can use the well-known exact formulae involving the Bernoulli numbers $B_k$:
\(
\zeta(2\ell) =  (-1)^{\ell+1} \frac{B_{2\ell}(2\pi)^{2\ell}} {2(2\ell)!},
\)
for every $\ell\in \N$, $\ell\ge 1$, where  the Bernoulli numbers $B_k$ are  defined
as the coefficients of the following series expansion:
\(
\frac{t}{e^t-1} =  \sum_{k=0}^{+\infty} B_k\frac{t^k}{k!},
\)
$\vert t \vert <2\pi$,
see, \emph{e.g.}, Cohen's book \cite[chapter 9]{Cohen2007}.
\end{enumerate}

 \subsection{Comparing running times}
 We implemented  \eqref{Gamma-x>1/2}-\eqref{Gamma-x<1/2}  and
 the formulae of Proposition \ref{DIF-formulae-gamma}  both in the
 scripting language of PARI/GP and in the C programming language. 
 In the first case (PARI/GP and gp2c),  using a precision of $128$ bits,
 \emph{i.e}, letting $n=128$, we compared
 the practical running times of computing $\log \Gamma (g^k/q)$, $k=0,\dotsc,q-2$,
 for $q=10007,305741,6766811,10000019,28227761$, $g$ being a fixed primitive root of $q$.
 In all these cases the use of \eqref{Gamma-x>1/2}-\eqref{Gamma-x<1/2} improved the total running
 times by a  40\%  factor with respect to the ones obtained using the predefined functions of PARI/GP. 
 Further improvements can be obtained using Proposition  \ref{DIF-formulae-gamma}
 if the particular application we are working on allows its use.  
 
In the second  case (C programming language), we repeated the computation previously described and 
we then compared the running times of 
our  implementation of \eqref{Gamma-x>1/2}-\eqref{Gamma-x<1/2}
and of the  \emph{long double} precision
version of $\log \Gamma$ defined in the C  language
(the \texttt{lgammal} function).
We clearly used a precision of $80$ bits
($n=80$). In this case our functions are slower of a factor $2.5$
than \texttt{lgammal}, while,
for the formulae of Proposition  \ref{DIF-formulae-gamma},
our functions are  slower of a factor $1.3$ with respect to \texttt{lgammal}. 
In both cases a  low-level  implementation of
our results might lead to a different outcome.

 \begin{Remark}[A digression on $\sin u$, $\pi$ and $\gamma$]
\label{sin-zeta}
 \begin{enumerate}[1)]
 \item
 A straightforward computation which uses $\Gamma(x) \Gamma(1-x)  = \pi / \sin(\pi x)$
 and Proposition \ref{DIF-formulae-gamma}  immediately gives 
 the well-known formula
 \begin{align}
\label{sin-x<1/2}
\sin(\pi x)
=
 \pi x   
\exp \Big(-\sum_{\ell=1}^{+\infty} \frac{\zeta(2\ell)}{\ell} x^{2\ell}
\Bigr) ,
\quad (0<x\le\frac12),
\end{align}
which, combined with the parity of the  $\sin$-function  and $\sin (0)=0$, can also 
be extended  to the whole interval $-1/2<x\le 1/2$. Equation \eqref{sin-x<1/2}, 
combined with the Bernoulli numbers definition, gives also
\[
\sin u 
=
u   \exp \Big(\sum_{\ell=1}^{+\infty} (-1)^{\ell}
 \frac{2^{2\ell-1}} {\ell(2\ell)!} B_{2\ell} u^{2\ell}
\Bigr) ,
\]
for every $u\in(-\pi/2,\pi/2]$.
    \item 
Computation of $\pi$. As a matter of curiosity, since we know that 
there are faster algorithms for this task, we remark that,
computing $\log\Gamma(x)+\log\Gamma(1-x)$ at $x=1/2$ with  
\eqref{Gamma-x>1/2}-\eqref{Gamma-x<1/2},
we obtain 
\[
\log \pi
= \log 2
+
\sum_{\ell=1}^{+\infty} \frac{\zeta(2\ell)}{\ell 4^\ell}   
\]
which can be used to compute $\pi$;  in fact, a straightforward implementation
using the scripting language of PARI/GP let us compute $10\,000$ decimal 
digits of $\pi$ in about $3$ seconds and $277$  milliseconds
while, for getting $1\,000$ decimal digits we just needed $19$ milliseconds
on the Dell OptiPlex machine previously mentioned.  
\item
 Computation of $\gamma$.
  As a matter of curiosity, since we know that 
 there are faster algorithms for this task, we remark that,
computing $\log\Gamma(x)-\log\Gamma(1-x)$ at $x=1/2$ with  
\eqref{Gamma-x>1/2}-\eqref{Gamma-x<1/2},
we obtain  the following  result  (first obtained by Stieltjes in 1887):
\[
\gamma
=  \log 2
-
\sum_{\ell=1}^{+\infty} \frac{\zeta(2\ell+1)}{(2\ell+1) 4^\ell} .  
\]
Such last formula can be clearly used to   compute $\gamma$; in fact, a 
straightforward implementation
using the scripting language of PARI/GP let us compute $10\,000$ decimal 
digits of $\gamma$ in about  $5$ minutes, $19$ seconds and  $255$ milliseconds
while, for getting $1\,000$ decimal digits we just needed $167$ milliseconds
on the Dell OptiPlex machine mentioned before.  
\end{enumerate}
 \end{Remark}

\section{An alternative algorithm to compute $\psi(x)$, $x\in(0,1)$}
\label{psi-computation}

Here we apply
to the \emph{digamma} function $\psi(x)=\Gamma^\prime/\Gamma(x)$, $x\in(0,1)$,
 the same argument used in Section \ref{loggamma-computation}. 
 The starting point
is the well-known Euler formula (see, \emph{e.g.}, Lagarias \cite[section 3]{Lagarias2013}):
\begin{equation}
\label{Euler-psi} 
\psi(x)  = - \gamma  - \sum_{k=2}^{+\infty} \zeta(k) (1-x)^{k-1} . 
\end{equation}

We immediately remark that the series in \eqref{Euler-psi}
absolutely converges for $x\in(0,2)$; this fact and the well-known relation 
\begin{equation}
\label{difference-psi}
\psi (1+x) = \psi (x)  +  \frac1x
\end{equation}
 let us 
obtain $\psi(x)$,  $x\in (0,1)$, in two different ways.
Recalling $\psi(1)=-\gamma$ and $\psi(1/2)=-2\log 2 -\gamma$,
we also remark that, letting $n\in \N$, $n\ge 2$, 
for every $x\in(0,2)$ there exists $r=r_\psi(x,n)  \ge 2$ 
such that
\begin{equation}
\label{tail-euler-series-psi}
\Bigl\vert \sum_{k=r+1}^{+\infty} \zeta(k) (1-x)^{k-1}  \Bigr\vert
<
 \zeta(3) 
\sum_{k=r+1}^{+\infty} 
\vert 1-x \vert^{k-1}
=
 1.21
\frac{\vert 1-x \vert^{r}}{1-\vert 1-x\vert} 
<2^{-n-1}.
\end{equation}
A straightforward computation reveals that we can choose
\begin{equation*}
r_\psi(x,n) = \Bigl\lceil  \frac{(n+2) \log 2 + \vert \log (1-\vert 1-x\vert)\vert}{ \vert \log  \vert 1-x \vert \vert}\Bigr\rceil.
\end{equation*}

\subsection{The shifting trick for $\psi(x)$, $x\in (0,1)$}
As for $r_\Gamma(x,n)$, we clearly have that $r_\psi(x,n)$ becomes larger as $\vert 1-x \vert$ increases.
We also remark that, using \eqref{difference-psi}, we can exploit the shifting trick in this case too. 
This way we will always use the best convergence interval, $x\in(1/2,3/2)$, we have for
the series in \eqref{Euler-loggamma}; we also remark that 
$r_\psi(x,n) \le r_\psi(1/2,n) =r_\psi(3/2,n) =  n+3$ for every $x\in(1/2,3/2)$.  
Summarising,  using  \eqref{Euler-psi} and \eqref{tail-euler-series-psi},
 for $x\in(1/2,1)$ we have that there exists $\theta=\theta(x) \in (-1/2,1/2)$ such that
 \begin{equation} 
 \label{psi-x>1/2}
\psi(x)  
= - \gamma  - \sum_{k=2}^{+\infty} \zeta(k) (1-x)^{k-1} 
= - \gamma  - \sum_{k=2}^{r_\psi(x,n)} \zeta(k) (1-x)^{k-1} 
 +\vert \theta \vert 2^{-n}. 
  \end{equation}
  We also  remark that for  $x\in(1/2,1)$, we have 
\(
 r_\psi(x,n) =\lceil \frac{(n+2) \log 2 + \vert \log x\vert}{ \vert \log  ( 1-x ) \vert} \rceil \le n+3.
\)
Moreover, using  \eqref{Euler-psi}-\eqref{tail-euler-series-psi}, 
for $x\in(0,1/2)$ we have  that there exists $\eta=\eta(x) \in (-1/2,1/2)$ such that
 \begin{equation}
 \label{psi-x<1/2}
\psi(x) 
=
- \frac{1}{x}  -\gamma   - \sum_{k=2}^{+\infty} (-1)^{k-1} \zeta(k) x^{k-1}
=
- \frac{1}{x}  -\gamma   - \sum_{k=2}^{r^\prime_\psi(x,n)} (-1)^{k-1} \zeta(k) x^{k-1}
+\vert \eta \vert 2^{-n},  
  \end{equation}
  where $r^\prime_\psi(x,n):=r_\psi(1+x,n) =  
  \lceil\frac{(n+2) \log 2 + \vert \log (1-  x )\vert}{ \vert \log    x   \vert}\rceil  \le n+3$.

We also remark that the series in the middle of  \eqref{psi-x>1/2}-\eqref{psi-x<1/2} can also
be obtained from the ones in  \eqref{Gamma-x>1/2}-\eqref{Gamma-x<1/2} by differentiation.
Since the needed $\zeta$-values  can be precomputed and stored with the desired precision
(using, for example, PARI/GP), 
 the formulae on the right hand sides of \eqref{psi-x>1/2}-\eqref{psi-x<1/2} allow us to compute
$\psi(x)$, $x\in (0,1)$, with a precision of $n$ binary digits 
using at most $n+3$ summands; moreover they also reveal that, from a computational point of view,
such a task is essentially as difficult
as computing   $1/ x$ when $x$ is close to $0$.

\begin{Remark}[Computation in the whole real axis]
It is clear that using \eqref{difference-psi} and \eqref{psi-x>1/2}-\eqref{psi-x<1/2}
we can compute $\psi(x)$ for every $x>0$ as follows.
For every $x>0$, we denote
as $\lfloor x \rfloor$ the integral part of $x$ and as $\{x\} = x - \lfloor x \rfloor$
the fractional part of $x$. Hence we obtain:
\begin{enumerate}[i)]
\item
 $\psi(1)=-\gamma$  and   $\psi(m)  = -\gamma +\sum_{k=1}^{m-1} 1/k$
for every $m\in \N$, $m\ge 2$;
\item
for  $x>1$, $x\not \in \N$, we  compute $\psi(x)$ as $\psi(x) = \psi(\{x\}) 
+ \sum_{k=0}^{\lfloor x \rfloor - 1}1/ (\{x\}+k)$;
\item
$\psi(1/2)=-2\log 2 -\gamma$;
\item
for $x\in (0,1/2)$,  we  compute $\psi(x)$ as in \eqref{psi-x<1/2};
\item
for  $x\in (1/2,1)$,   we  compute $\psi(x)$ as in \eqref{psi-x>1/2}.
\end{enumerate}
Even if we are mainly interested in working with $x\in(0,1)$ we recall that for $x$ large it 
might be more efficient to implement an asymptotic formula for $\psi(x)$.
\end{Remark}

\begin{Remark} [Enlarging the convergence radius] The size of the convergence interval  in \eqref{Euler-psi} 
can be doubled by isolating the Taylor series at $1$ of $1/x -1$ in \eqref{Euler-psi} 
thus getting
\begin{equation}
\label{enlarged-radius-psi}
\psi(x) = 
-\frac{1}{x}- \gamma + 1 - \sum_{k=2}^{+\infty} (\zeta(k)-1) (1-x)^{k-1} .
\end{equation}
Using the well-known estimate
$\vert \zeta(k)-1 \vert<2^{1-k}$ for every $k\in\N, k\ge 3$,
it is easy to prove that  the  series in \eqref{enlarged-radius-psi} converges for every $x\in(-1,3)$.  
\end{Remark}

\begin{Remark}[Computation in the complex plane]
\begin{enumerate}[1)]
\item
It seems that the argument leading to \eqref{psi-x>1/2}-\eqref{psi-x<1/2}
is not usually implemented in the most used software libraries or Computer 
Aided Systems for Mathematics for the same reasons we discussed for the $\log\Gamma$-function. 
\item
We argue analogously as we did for the $\log \Gamma$-function.
Using the formula 
\begin{equation}
\label{psi-complex}
\psi(z) 
= 
- \gamma
  - \sum_{k=2}^{+\infty} \zeta(k)(1-z)^{k-1} 
\end{equation}
which generalises \eqref{Euler-psi} to the region $\vert z-1\vert<1$,
and \eqref{difference-psi}, which in fact holds for any argument 
$z\in \C$, $z\ne -n$, $n\in \N$, we can build a possible complex strip
 starting from the rectangle 
$1/2\le \Re(z)< 3/2$, $\vert \Im(z)\vert \le 1/4$,
in which less than 
$1.2\cdot(n+2)+4$ terms
are sufficient to have a precision of $n$ binary digits in computing
a suitable truncation of the  series
in \eqref{psi-complex}. 
Moreover, such a strip can be  vertically
enlarged using the following form of Gau\ss' multiplication
theorem 
\[
\psi(mz) =
\log m
+ 
\frac{1}{m}
\sum_{j=0}^{m-1}\psi\Bigl(z+\frac{j}{m}\Bigr), 
\]
where $m\in \N$, $m\ge 1$.
\end{enumerate}
\end{Remark}

\subsection{Reflection formulae}
We now remark that in
using    \eqref{psi-x>1/2}-\eqref{psi-x<1/2} to compute $\psi(x) -\psi(1-x)$,
$x\in (0,1)$,   the odd summands in the series will vanish
(and in the corresponding series for $\psi(x) + \psi(1-x)$ the even summands will be discarded too).
We summarise the situation in the following
\begin{Proposition}
\label{DIF-formulae-psi}
Let $x\in (0,1)$, $x\ne 1/2$,
$n\in \N$, $n\ge 2$, $r_1(x,n)=\lceil\frac{(n+2) \log 2 + \vert \log (1-x) \vert}{ \vert \log  x \vert }\rceil/2$
and $r_2(x,n)  = \lceil\frac{(n+2) \log 2 +\vert \log x \vert}{\vert \log  (1-x) \vert}\rceil/2$.
Using  \eqref{psi-x>1/2}-\eqref{psi-x<1/2},
we have that there exists $\theta=\theta(x)\in (-1/2,1/2)$  such that
for $0<x <1/2$ we have
\begin{align}  
\psi(x)&-\psi(1-x) 
 \label{psi-DIF-x<1/2}  
=
- \frac{1}{x}  +
2\sum_{\ell=1}^{r_1} \zeta(2\ell) x^{2\ell-1} 
+ \vert \theta \vert 2^{-n} ,
\\
\psi(x)&+\psi(1-x) 
 \label{psi-DIF-even-x<1/2}  
=
-2\gamma
-\frac{1}{x}  - 
2\sum_{\ell=1}^{r_1} \zeta(2\ell+1) x^{2\ell} 
+ \vert \theta \vert 2^{-n},
\end{align}
and for $1/2<x <1$ we have
  \begin{align}
\psi(x)&-\psi(1-x) 
\label{psi-DIF-x>1/2} 
=
\frac{1}{1-x}  -
2\sum_{\ell=1}^{r_2}  \zeta(2\ell) (1-x)^{2\ell-1} 
+ \vert \theta \vert 2^{-n} ,  
\\ 
\psi(x)&+\psi(1-x) 
\label{psi-DIF-even-x>1/2} 
=
-2\gamma
-\frac{1}{1-x}  
-
2\sum_{\ell=1}^{r_2}  \zeta(2\ell+1) (1-x)^{2\ell} 
+ \vert \theta \vert 2^{-n} .
 \end{align}
\end{Proposition}
 \begin{Proof}
Assume that $0<x < 1/2$; we compute $\psi(x)$ with the series in  \eqref{psi-x<1/2} 
and $\psi(1-x)$ with the series in   \eqref{psi-x>1/2}.
Since they absolutely converge, their sum is obtained with the series
having as summands the sum of their coefficients. Arguing as in 
\eqref{tail-euler-series-psi} and remarking that  
$r_1(x,n)= r_\psi(1-x,n) /2=  r^\prime_\psi(x,n)/2$,
we immediately have  that \eqref{psi-DIF-x<1/2} holds since
  the odd summands vanish.
Assume that $1/2<x<1$; in this case we  compute $\psi(x)$ with the series in  \eqref{psi-x>1/2} 
and $\psi(1-x)$ with the series in   \eqref{psi-x<1/2}. 
Since they absolutely converge, their sum is obtained with the series
having as summands the sum of their coefficients. Arguing as in 
\eqref{tail-euler-series-psi} and remarking that  
$r_2(x,n)= r_\psi(x,n)/2 =  r^\prime_\psi(1-x,n)/2$, we
immediately have  that \eqref{psi-DIF-x>1/2} holds since
  the odd summands vanish.  
  The  derivation of \eqref{psi-DIF-even-x<1/2} and \eqref{psi-DIF-even-x>1/2}
 is similar.
 This completes the proof.
\end{Proof}
 
It is worth mentioning that the right hand side in \eqref{psi-DIF-x>1/2} can
be obtained from the one in \eqref{psi-DIF-x<1/2} formally replacing $x$ with $1-x$
and changing sign;
and that, without the change of sign, the same holds for  \eqref{psi-DIF-even-x<1/2}    
and \eqref{psi-DIF-even-x>1/2}.
Using $\psi(1-x) -\psi(x) = \pi  \cot(\pi x)$,   Proposition \ref{DIF-formulae-psi}
 also immediately gives a way of writing $ \cot(\pi x)$
in term of logs and values of the Riemann zeta-function at positive even integers,
see Remark \ref{cot-zeta} below. 
The use of Proposition \ref{DIF-formulae-psi} 
is particularly efficient for 
the same reasons we already described for Proposition \ref{DIF-formulae-gamma};
we just need to remark that the number of summands in this 
case is  $\le (n+3)/2$.

 \begin{Remark}[A digression on $\cot u$]
 \label{cot-zeta}
 A straightforward computation which uses $\psi(1-x) -\psi(x) = \pi  \cot(\pi x)$
 and Proposition \ref{DIF-formulae-psi}  immediately gives 
 the well-known formula
\[ 
\cot(\pi x)
=
\frac{1}{\pi x}  - 
\frac{2}{\pi}
\sum_{\ell=1}^{+\infty} \zeta(2\ell) x^{2\ell-1}  
\quad (0<x\le\frac12),
\]
which, combined with the parity of the cotangent function, can also be extended 
to $-1/2<x\le 1/2$, $x\ne 0$.
 \end{Remark}

\renewcommand{\bibliofont}{\normalsize} 

\vskip 1cm
\noindent
Alessandro Languasco
Universit\`a di Padova,
 Dipartimento di Matematica,
 ``Tullio Levi-Civita'',
Via Trieste 63,
35121 Padova, Italy.   
{\it e-mail}: alessandro.languasco@unipd.it \\

\newpage 

\begin{table}[htp]
\scalebox{0.606}{
\begin{tabular}{|c|c|}
\hline
$q$  &  $M_q$\\ \hline 
 $3$ & $0.60459978807807261686469275254\dotsc$ \\
$5$ & $0.88857658763167324940317619801\dotsc$ \\
$7$ & $1.18741041172372594878462529795\dotsc$ \\
$11$ & $1.42640418224108352050983157388\dotsc$ \\
$13$ & $1.40613477980703732992641904009\dotsc$ \\
$17$ & $1.64849370699838393605712538480\dotsc$ \\
$19$ & $1.66331503401599345646761861198\dotsc$ \\
$23$ & $1.96520205410785916590276700512\dotsc$ \\
$29$ & $1.94760835256298812067787472221\dotsc$ \\
$31$ & $1.95460789577555396208826644691\dotsc$ \\
$37$ & $1.99550481523309639952094476513\dotsc$ \\
$41$ & $2.11431182971789691740711719452\dotsc$ \\
$43$ & $2.15300367350566872872181166589\dotsc$ \\
$47$ & $2.29124192852861593669991478644\dotsc$ \\
$53$ & $2.30607194293623581454960238070\dotsc$ \\
$59$ & $2.37483432765382947109503698217\dotsc$ \\
$61$ & $2.27383316907813941451707465286\dotsc$ \\
$67$ & $2.40177874951129444901523775053\dotsc$ \\
$71$ & $2.60986917715784586434512887899\dotsc$ \\
$73$ & $2.22099352575696724500691589494\dotsc$ \\
$79$ & $2.49865862902662621045751915333\dotsc$ \\
$83$ & $2.40279523907172221214735958325\dotsc$ \\
$89$ & $2.48752834330666367191078756009\dotsc$ \\
$97$ & $2.42085614235917869433795064647\dotsc$ \\
$101$ & $2.49348309598992905601857403054\dotsc$ \\
$103$ & $2.58872219777793564220546520282\dotsc$ \\
$107$ & $2.54845309686851323582504835532\dotsc$ \\
$109$ & $2.43977808110771365276355432077\dotsc$ \\
$113$ & $2.37003547987807971651823457293\dotsc$ \\
$127$ & $2.72368730766675765849410753531\dotsc$ \\
$131$ & $2.57848536120487307251550163471\dotsc$ \\
$137$ & $2.72051117298997095021332225039\dotsc$ \\
$139$ & $2.78392266035840572927430273739\dotsc$ \\
$149$ & $2.68631183937556722239373025709\dotsc$ \\
$151$ & $2.62352855587439821150270052709\dotsc$ \\
$157$ & $2.91562362895517732081839286279\dotsc$ \\
$163$ & $2.69099736683125370993373921388\dotsc$ \\
$167$ & $2.74644085264695532114443109623\dotsc$ \\
$173$ & $2.83083393377236324187600452713\dotsc$ \\
$179$ & $2.95215347090837063371989208638\dotsc$ \\
$181$ & $2.55866549759635341623707612755\dotsc$ \\
$191$ & $2.95512966360404799352636309788\dotsc$ \\
$193$ & $2.60255291569166233786515416685\dotsc$ \\
$197$ & $2.81468933588096728324501080140\dotsc$ \\
$199$ & $2.79249308566493928174043020396\dotsc$ \\
$211$ & $2.89594376660833329272394923676\dotsc$ \\
$223$ & $2.92740993747063127172267612416\dotsc$ \\
$227$ & $2.68264264675366168697033571053\dotsc$ \\
$229$ & $2.87581890996867882018490166417\dotsc$ \\
$233$ & $2.91527084775689862127414734394\dotsc$ \\
$239$ & $3.04819103378239805449031098106\dotsc$ \\
$241$ & $2.77830595707238913568405847443\dotsc$ \\
$251$ & $2.96609382990202814045159294394\dotsc$ \\
$257$ & $2.90271693301413614997407956341\dotsc$ \\
$263$ & $2.93609043858561342569242914697\dotsc$ \\
$269$ & $2.95235925085763617694817528815\dotsc$ \\
\hline
\end{tabular}
}
\scalebox{0.606}{
\begin{tabular}{|c|c|}
\hline
$q$  &  $M_q$\\ \hline
$271$ & $2.86925767656885353188813607888\dotsc$ \\
$277$ & $3.01622855659087525191859838663\dotsc$ \\
$281$ & $2.80321128554057954903648560934\dotsc$ \\
$283$ & $2.94572301912315375394664649007\dotsc$ \\
$293$ & $3.11291787542876229268348893252\dotsc$ \\
$307$ & $3.02482538363473232115657507118\dotsc$ \\
$311$ & $3.38472414241331604469710895130\dotsc$ \\
$313$ & $3.10929623072433076283321530583\dotsc$ \\
$317$ & $2.95307987036615177805014754274\dotsc$ \\
$331$ & $2.95426729345152158301122222583\dotsc$ \\
$337$ & $3.01307766294270421990892299882\dotsc$ \\
$347$ & $3.10750424169246830621098007832\dotsc$ \\
$349$ & $3.24312632555570194704840092849\dotsc$ \\
$353$ & $2.96318053008948669143319019321\dotsc$ \\
$359$ & $3.15033145392974552469052671095\dotsc$ \\
$367$ & $2.94003500276564543524572051749\dotsc$ \\
$373$ & $3.32946093758130126984874928456\dotsc$ \\
$379$ & $2.98417634686154821282891886779\dotsc$ \\
$383$ & $3.15698301880900757107584482170\dotsc$ \\
$389$ & $3.13576641023433071248884490813\dotsc$ \\
$397$ & $3.22908440562028461328032006458\dotsc$ \\
$401$ & $3.02531335084626851927421904169\dotsc$ \\
$409$ & $2.97468553326599915860870900047\dotsc$ \\
$419$ & $3.15638028965444731073868407588\dotsc$ \\
$421$ & $2.89771782676987997033822961589\dotsc$ \\
$431$ & $3.27666788721708595129425883251\dotsc$ \\
$433$ & $3.14357057349630992081515342345\dotsc$ \\
$439$ & $3.22131456032260084903381293645\dotsc$ \\
$443$ & $3.28708635670026587878865275729\dotsc$ \\
$449$ & $3.14374861185339258666718155942\dotsc$ \\
$457$ & $3.01805880291502467267265047463\dotsc$ \\
$461$ & $2.99061762375575687939818599075\dotsc$ \\
$463$ & $3.02557292452165996572378279291\dotsc$ \\
$467$ & $3.17549334417045753275838551244\dotsc$ \\
$479$ & $3.58857580472017716180065757364\dotsc$ \\
$487$ & $3.02150659065323210851247532813\dotsc$ \\
$491$ & $3.08880305445418175347931846454\dotsc$ \\
$499$ & $3.33310237076133842063983002611\dotsc$ \\
$503$ & $3.23862981619088362793681328519\dotsc$ \\
$509$ & $3.18790240372285890369329278804\dotsc$ \\
$521$ & $3.32391248073785183118435027656\dotsc$ \\
$523$ & $3.39266087785391065974002131993\dotsc$ \\
$541$ & $3.08990995669354475231608583966\dotsc$ \\
$547$ & $3.21177134746313664845285051418\dotsc$ \\
$557$ & $3.36888198211914580127389042801\dotsc$ \\
$563$ & $3.19827255926659391012489648438\dotsc$ \\
$569$ & $3.19231644611870880739554354757\dotsc$ \\
$571$ & $3.29737963928860833400612713236\dotsc$ \\
$577$ & $3.10129116198243754093122394033\dotsc$ \\
$587$ & $3.17339593032303497730382120451\dotsc$ \\
$593$ & $3.32734104110339616101985098272\dotsc$ \\
$599$ & $3.20904989883664076353720703220\dotsc$ \\
$601$ & $3.03462372828981497636351602452\dotsc$ \\
$607$ & $3.22334095612061763642992301825\dotsc$ \\
$613$ & $3.30461730437919519650689110196\dotsc$ \\
$617$ & $3.21110531139170629353622676966\dotsc$ \\
\hline
\end{tabular}
}
\scalebox{0.606}{
\begin{tabular}{|c|c|}
\hline
$q$  &  $M_q$\\ \hline
  $619$ & $3.37069356094491605525491036082\dotsc$ \\
$631$ & $3.13208638787415093755171585220\dotsc$ \\
$641$ & $3.23211518539269050934542501919\dotsc$ \\
$643$ & $3.18685340027257900687187281884\dotsc$ \\
$647$ & $3.38314148149622511105230281091\dotsc$ \\
$653$ & $3.43910390260448591467634183370\dotsc$ \\
$659$ & $3.40474114893667565054038830288\dotsc$ \\
$661$ & $3.19707359599057009574455632193\dotsc$ \\
$673$ & $3.29524177279177739366473959234\dotsc$ \\
$677$ & $3.14203413621204495006097125621\dotsc$ \\
$683$ & $3.10413178955563308313872650724\dotsc$ \\
$691$ & $3.51698989024002615984946445256\dotsc$ \\
$701$ & $3.35085324504731590551613080436\dotsc$ \\
$709$ & $3.18491137055861253375845988616\dotsc$ \\
$719$ & $3.63201071708524524489955813645\dotsc$ \\
$727$ & $3.19488076113527166415118327181\dotsc$ \\
$733$ & $3.45685741198240493832871170689\dotsc$ \\
$739$ & $3.29109339411323807636519499288\dotsc$ \\
$743$ & $3.36194027636618645067801019522\dotsc$ \\
$751$ & $3.37121298850114852069195462707\dotsc$ \\
$757$ & $3.29522003757629360256760087960\dotsc$ \\
$761$ & $3.39036396026169168862024302997\dotsc$ \\
$769$ & $3.37966794694746138163416821630\dotsc$ \\
$773$ & $3.32182777438635162483294645420\dotsc$ \\
$787$ & $3.41806860456589489290488883781\dotsc$ \\
$797$ & $3.41332527420052976238715906739\dotsc$ \\
$809$ & $3.33399463072530149065220868852\dotsc$ \\
$811$ & $3.41748712961750993035612332531\dotsc$ \\
$821$ & $3.32756221988730824824951614294\dotsc$ \\
$823$ & $3.37960627332789369390005691334\dotsc$ \\
$827$ & $3.40849620130058074419186199017\dotsc$ \\
$829$ & $3.37089500117347869458163604211\dotsc$ \\
$839$ & $3.57917416660357193379948511931\dotsc$ \\
$853$ & $3.38016417685190138464013883123\dotsc$ \\
$857$ & $3.57544465806174330360312716184\dotsc$ \\
$859$ & $3.50541693805340088235348934788\dotsc$ \\
$863$ & $3.50825623636559589065557628386\dotsc$ \\
$877$ & $3.49769275741261662405276790520\dotsc$ \\
$881$ & $3.39897242344112875811946433945\dotsc$ \\
$883$ & $3.34652632278953040018003375892\dotsc$ \\
$887$ & $3.34968741972912759452502545654\dotsc$ \\
$907$ & $3.39778872837837091734236419151\dotsc$ \\
$911$ & $3.37032660476473849985853056954\dotsc$ \\
$919$ & $3.35063473135683107552976812180\dotsc$ \\
$929$ & $3.53321177004124872239874436148\dotsc$ \\
$937$ & $3.56404261890133745461623907397\dotsc$ \\
$941$ & $3.36178737964711744410945269493\dotsc$ \\
$947$ & $3.35241151685426146046317817853\dotsc$ \\
$953$ & $3.49561078035208420652588273363\dotsc$ \\
$967$ & $3.34475180524849073411075928718\dotsc$ \\
$971$ & $3.50133666092278415387426395493\dotsc$ \\
$977$ & $3.46323541639460242041962117061\dotsc$ \\
$983$ & $3.38200836946693885301545542710\dotsc$ \\
$991$ & $3.47683493471726594626831378059\dotsc$ \\
$997$ & $3.49595691818271364657373564526\dotsc$ \\
  \phantom{} & \phantom{}  \\
\hline
\end{tabular}
}
\caption{\label{table1}
{\small
Values of $M_q$ for every odd prime up to $1000$ with 
$30$-digit precision; computed with PARI/GP, v.~2.11.4,
with  trivial summing over $a$.
Total computation time:  16 sec., 539 millisecs. 
}
}
\end{table} 
\newpage

\begin{table}[htp]
\scalebox{0.575}{
\begin{tabular}{|c|c|}
\hline
$q$  &  $m_q$\\ \hline 
 $3$ & $0.604599788078072616864692752547\dotsc$ \\
$5$ & $0.430408940964004038889433232951\dotsc$ \\
$7$ & $0.547959686797993973084485988763\dotsc$ \\
$11$ & $0.618351934876807874060419662662\dotsc$ \\
$13$ & $0.598987497945465758207499250242\dotsc$ \\
$17$ & $0.453546340908733659287599108349\dotsc$ \\
$19$ & $0.413193436540565451244291268589\dotsc$ \\
$23$ & $0.552304916713058385866569568830\dotsc$ \\
$29$ & $0.451093787499735920935796126723\dotsc$ \\
$31$ & $0.440122223433962808617040155019\dotsc$ \\
$37$ & $0.420012687371836528710987165573\dotsc$ \\
$41$ & $0.531058786094205975539953409040\dotsc$ \\
$43$ & $0.479088388239857211764493892920\dotsc$ \\
$47$ & $0.367129807516487530311860512131\dotsc$ \\
$53$ & $0.413967522107356708589184273274\dotsc$ \\
$59$ & $0.332420580251333196195200170353\dotsc$ \\
$61$ & $0.365767556044524607545195023959\dotsc$ \\
$67$ & $0.383806628882915516388036615854\dotsc$ \\
$71$ & $0.445439715279504151396753278541\dotsc$ \\
$73$ & $0.332816572422086433683238786444\dotsc$ \\
$79$ & $0.428163851610403317805757592040\dotsc$ \\
$83$ & $0.462387538865549865563618907321\dotsc$ \\
$89$ & $0.389534102336872450885689912904\dotsc$ \\
$97$ & $0.389092591237496211100095150861\dotsc$ \\
$101$ & $0.356508197108894602900482691212\dotsc$ \\
$103$ & $0.397853727338476109339913228255\dotsc$ \\
$107$ & $0.404936054291382959432040324919\dotsc$ \\
$109$ & $0.395439291752596225073889921584\dotsc$ \\
$113$ & $0.429318876983891967381982481536\dotsc$ \\
$127$ & $0.370583958563266768017674320198\dotsc$ \\
$131$ & $0.352086298168602566529507073368\dotsc$ \\
$137$ & $0.417505824794320762685433840755\dotsc$ \\
$139$ & $0.416862458745916605478311229554\dotsc$ \\
$149$ & $0.383901386511619291505439140039\dotsc$ \\
$151$ & $0.375982979158935926530947639675\dotsc$ \\
$157$ & $0.369949970311229129154525082455\dotsc$ \\
$163$ & $0.246068527552960243897853273760\dotsc$ \\
$167$ & $0.350631724723697517489360519141\dotsc$ \\
$173$ & $0.382663127794428634867406502671\dotsc$ \\
$179$ & $0.341096922920489890524209012663\dotsc$ \\
$181$ & $0.400130740444375343545002983965\dotsc$ \\
$191$ & $0.369946585080866896719707254253\dotsc$ \\
$193$ & $0.362378327133943195753957944776\dotsc$ \\
$197$ & $0.373837653864628774889214137854\dotsc$ \\
$199$ & $0.407431952492485762271173360682\dotsc$ \\
$211$ & $0.352658667943486049139179190376\dotsc$ \\
$223$ & $0.430673753410189075458236898964\dotsc$ \\
$227$ & $0.381439782634825445136263631484\dotsc$ \\
$229$ & $0.383729157778237776053663569946\dotsc$ \\
$233$ & $0.365137320866318953249895485603\dotsc$ \\
$239$ & $0.361752701714412120902910350331\dotsc$ \\
$241$ & $0.352065536871065140512991463133\dotsc$ \\
$251$ & $0.423961989644231760650605953707\dotsc$ \\
$257$ & $0.388887492169207487425676370870\dotsc$ \\
$263$ & $0.348202135404538171971904967311\dotsc$ \\
$269$ & $0.404844287635172668787118855557\dotsc$ \\
\hline
\end{tabular}
}
\scalebox{0.575}{
\begin{tabular}{|c|c|}
\hline
$q$  &  $m_q$\\ \hline
 $271$ & $0.381676325868210243355238004400\dotsc$ \\
$277$ & $0.367772078360884977669216925848\dotsc$ \\
$281$ & $0.367654214871298109838839266841\dotsc$ \\
$283$ & $0.360060056162870590740646688749\dotsc$ \\
$293$ & $0.331438394291933793807551984375\dotsc$ \\
$307$ & $0.338116784161146375199978675010\dotsc$ \\
$311$ & $0.364786872431841488117413342623\dotsc$ \\
$313$ & $0.348283081331744875759624637812\dotsc$ \\
$317$ & $0.353920323300651284267391675261\dotsc$ \\
$331$ & $0.384662955943923874879364794182\dotsc$ \\
$337$ & $0.395558761496812443683554845556\dotsc$ \\
$347$ & $0.377945303983869233549712938160\dotsc$ \\
$349$ & $0.346392157694090105389329641835\dotsc$ \\
$353$ & $0.392050917804241245253197847943\dotsc$ \\
$359$ & $0.372107247747624127090483414376\dotsc$ \\
$367$ & $0.292520324658371008714768864693\dotsc$ \\
$373$ & $0.365107572487765521944814562236\dotsc$ \\
$379$ & $0.344529506663803122648629099168\dotsc$ \\
$383$ & $0.369414454892205779888789721588\dotsc$ \\
$389$ & $0.309685712040159744098308002254\dotsc$ \\
$397$ & $0.345438689174522431354681542746\dotsc$ \\
$401$ & $0.374008182395715337562657042856\dotsc$ \\
$409$ & $0.343912321305157763624260753609\dotsc$ \\
$419$ & $0.398860743207736694380625064238\dotsc$ \\
$421$ & $0.355558244246456895812338230951\dotsc$ \\
$431$ & $0.304249924284011002888004890483\dotsc$ \\
$433$ & $0.374437191561992583380814554550\dotsc$ \\
$439$ & $0.385201582268506040100915711884\dotsc$ \\
$443$ & $0.348072727822271656582475517939\dotsc$ \\
$449$ & $0.387387308651078457138994462088\dotsc$ \\
$457$ & $0.344128199358590128646593703075\dotsc$ \\
$461$ & $0.342282439127743493973555398134\dotsc$ \\
$463$ & $0.322810266925933801000377491409\dotsc$ \\
$467$ & $0.346262914614287495219714965304\dotsc$ \\
$479$ & $0.358217517993363694169353276634\dotsc$ \\
$487$ & $0.347056625920168854172557770589\dotsc$ \\
$491$ & $0.328216811995353442798930033747\dotsc$ \\
$499$ & $0.343916679837668925712980759633\dotsc$ \\
$503$ & $0.305270080648936874982633449014\dotsc$ \\
$509$ & $0.358400496109432533323994071825\dotsc$ \\
$521$ & $0.326639577620790076266875239423\dotsc$ \\
$523$ & $0.370263699010744524473643538398\dotsc$ \\
$541$ & $0.371841857227428197484414254772\dotsc$ \\
$547$ & $0.344451554730668674041547075694\dotsc$ \\
$557$ & $0.331277738629462528622134186640\dotsc$ \\
$563$ & $0.336288953230595448646092763257\dotsc$ \\
$569$ & $0.355458621282047158720345039619\dotsc$ \\
$571$ & $0.308093921260701873971115764321\dotsc$ \\
$577$ & $0.351768032853537169914900076369\dotsc$ \\
$587$ & $0.364218277748308523039650017021\dotsc$ \\
$593$ & $0.339509407115695414229878807307\dotsc$ \\
$599$ & $0.352374411311513496889433605686\dotsc$ \\
$601$ & $0.350992936007857504857261176829\dotsc$ \\
$607$ & $0.329981981645739789066558958285\dotsc$ \\
$613$ & $0.323436732518113096063105608197\dotsc$ \\
$617$ & $0.333733731209755325952269458510\dotsc$ \\
\hline
\end{tabular}
}
\scalebox{0.575}{
\begin{tabular}{|c|c|}
\hline
$q$  &  $m_q$\\ \hline 
$619$ & $0.363853513253924249520317635224\dotsc$ \\
$631$ & $0.325338069145847294024656423120\dotsc$ \\
$641$ & $0.352147403815645998874657459363\dotsc$ \\
$643$ & $0.337445372236761198609110503052\dotsc$ \\
$647$ & $0.361173357971713913910561639944\dotsc$ \\
$653$ & $0.344528546854302142286423852753\dotsc$ \\
$659$ & $0.338132160088261589245693055585\dotsc$ \\
$661$ & $0.360619845188122885583816505945\dotsc$ \\
$673$ & $0.376129394432057443687895596361\dotsc$ \\
$677$ & $0.303745675684779321753979753781\dotsc$ \\
$683$ & $0.333368382358596313476605760567\dotsc$ \\
$691$ & $0.383305396571474295606978030935\dotsc$ \\
$701$ & $0.357535298673103247072275901603\dotsc$ \\
$709$ & $0.358127149918438747730782304062\dotsc$ \\
$719$ & $0.370805494541737091847758990324\dotsc$ \\
$727$ & $0.372898383226526400669880140310\dotsc$ \\
$733$ & $0.350732615201881873478746719218\dotsc$ \\
$739$ & $0.304897212128635209913817219538\dotsc$ \\
$743$ & $0.352336063056891753760152372944\dotsc$ \\
$751$ & $0.328140563247974273370269669938\dotsc$ \\
$757$ & $0.319411975837870470401241096947\dotsc$ \\
$761$ & $0.321973143068848452624799526258\dotsc$ \\
$769$ & $0.338178333435205826998970942276\dotsc$ \\
$773$ & $0.354965105866232352811142020422\dotsc$ \\
$787$ & $0.324266393608051759348322750934\dotsc$ \\
$797$ & $0.350515221653708975482852748559\dotsc$ \\
$809$ & $0.326434531678769569498090895809\dotsc$ \\
$811$ & $0.308510351209093190875996531501\dotsc$ \\
$821$ & $0.330244313860677980685649524053\dotsc$ \\
$823$ & $0.324283030566918000839563175873\dotsc$ \\
$827$ & $0.329051480295616938819583708993\dotsc$ \\
$829$ & $0.321222221357102394251381440581\dotsc$ \\
$839$ & $0.351875438185464119441010492489\dotsc$ \\
$853$ & $0.351818805199888913965741566340\dotsc$ \\
$857$ & $0.366066564937826743322512006592\dotsc$ \\
$859$ & $0.331002403214058452880530388161\dotsc$ \\
$863$ & $0.307797465512278850449583310784\dotsc$ \\
$877$ & $0.340810810909356878124676806818\dotsc$ \\
$881$ & $0.333426901946204902333090622005\dotsc$ \\
$883$ & $0.317169031194076637573267354539\dotsc$ \\
$887$ & $0.338488672510289686224158744714\dotsc$ \\
$907$ & $0.312944615763902492538994990004\dotsc$ \\
$911$ & $0.333229470950277517416975539246\dotsc$ \\
$919$ & $0.355302031868156498197335032183\dotsc$ \\
$929$ & $0.357591086944711096303847153788\dotsc$ \\
$937$ & $0.340988515678272981620106866387\dotsc$ \\
$941$ & $0.340561178701543781023570327465\dotsc$ \\
$947$ & $0.379684850405948893784190778135\dotsc$ \\
$953$ & $0.341656471737723220834223917631\dotsc$ \\
$967$ & $0.320793935250234540306919187736\dotsc$ \\
$971$ & $0.305685071173978402484536717845\dotsc$ \\
$977$ & $0.342232473746414199291917612827\dotsc$ \\
$983$ & $0.346431365398324524472675676012\dotsc$ \\
$991$ & $0.343527153645248056548459399539\dotsc$ \\
$997$ & $0.307684966461541454308715912751\dotsc$ \\
 \phantom{} & \phantom{}  \\ 
\hline
\end{tabular}
}
\caption{\label{table2}
{\small
Values of $m_q$ for every odd prime up to $1000$ with 
$30$-digit precision; computed with PARI/GP, v.~2.11.4,
with  trivial summing over $a$.
Total computation time:  16 sec., 817 millisecs.  
}
}
\end{table}  
 
   \vfill\eject
\begin{figure}[ht] 
\includegraphics[scale=0.48,angle=0]{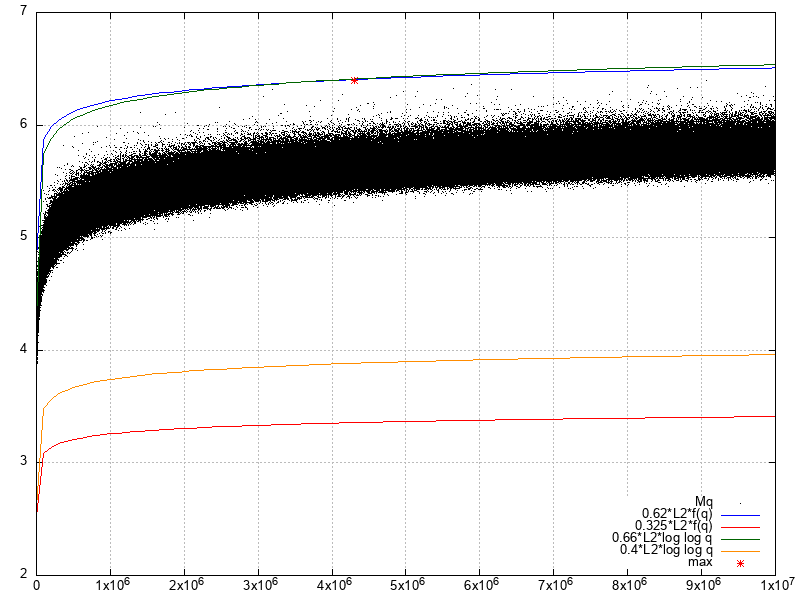}  
\caption{\small{Graphical representation of Theorem \ref{Thm-max}.
The values of $M_q$, $q$ prime, $500\le q\le \bound$.   
The  \captionmax\
The blue  line  represents $ 0.62 L_2  f(q)$;  the red one  
$0.325 L_2  f(q)$, where   $f(q)$ is defined in \eqref{fg-def}
and $L_2:=  2 e^\gamma$.
The green line represents $0.66L_2\log \log q$; the orange one $0.4L_2\log \log q$.
For $3\le q <500$, see the next plot. 
}
 }
\label{fig-LLS-1} 
 \end{figure}
 
\vfill\eject 

\vfill\eject 
 
 \vfill\eject 
\begin{figure}[ht] 
\includegraphics[scale=0.48,angle=0]{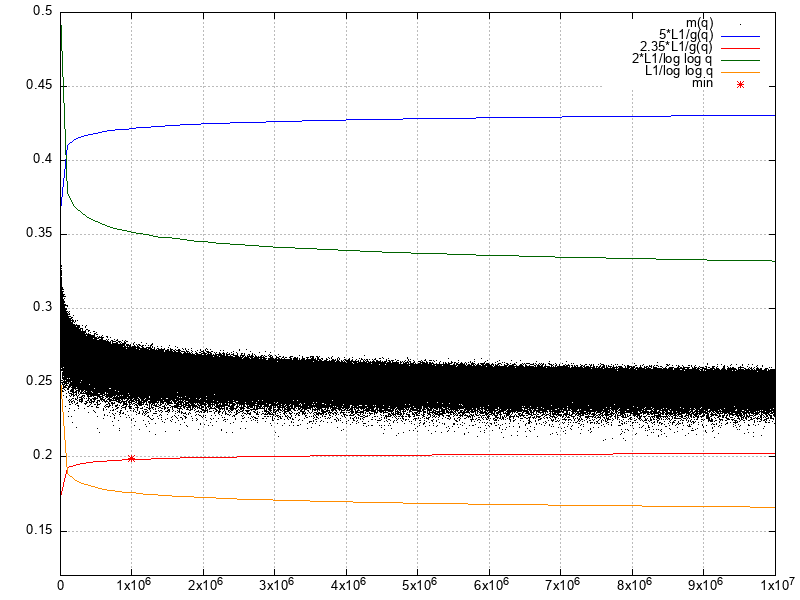}  
\caption{\small{Graphical representation of Theorem \ref{Thm-min}.
The values of $m_q$, $q$ prime, $500\le q\le \bound$.  
\captionmin\
The blue  line  represents $ 5 L_1 / g(q)$;  the red one  
$2.35 L_1 / g(q)$, where $g(q)$ is defined in \eqref{fg-def}
and $L_1:=   \frac{\pi^2}{12 e^\gamma}$. 
The green line represents $2L_1/\log\log q$; the orange one $1.13 L_1/\log\log q$.
 For $3\le q <500$, see the next plot.
}
 }
\label{fig-LLS-3} 
\end{figure} 

\vfill\eject 
 
\vfill\eject 
 
\end{document}